\documentclass[11pt]{article}
\usepackage{amsfonts, latexsym, amssymb,amsmath,amstext}
\usepackage[dvips]{graphics,color}
\usepackage[mathscr]{eucal}
\newenvironment{proof}{\par\noindent{\bf Proof.}\ }{\hfill$\Box$\par\medskip}
\newtheorem{theorem}{Theorem}[section]
\newtheorem{lemma}[theorem]{Lemma}

\newtheorem{example}[theorem]{Example}

\newcommand{\equal}{&\!\!\!=\!\!\!&}

\newcommand{\CC}{\mathbb C}
\newcommand{\PP}{\mathbb P}
\newcommand{\QQ}{\mathbb Q}
\makeatletter
\@addtoreset{equation}{section}

\makeatother
\begin{document}
\title{Quadratic Transformations of the Sixth Painlev\'e Equation
with Application to Algebraic Solutions}
\author{Raimundas Vid\= unas\thanks{Supported by the 21 Century COE Programme
"Development of Dynamic Mathematics with High Functionality" of the Ministry
of Education, Culture, Sports, Science and Technology of Japan. E-mail:
vidunas@math.kyushu-u.ac.jp}\; and\; Alexander~V.~Kitaev\thanks{Supported by
JSPS grant-in-aide no.~$14204012$. E-mail:
kitaev@pdmi.ras.ru}\\
Department of Mathematics, Kyushu University, 812-8581 Fukuoka, Japan\footnotemark[1]\\
Steklov Mathematical Institute, Fontanka 27, St. Petersburg, 191023, Russia\footnotemark[2]\\
and\\
School of Mathematics and Statistics, University of Sydney,\\
Sydney, NSW 2006, Australia\footnotemark[1]\;\,\footnotemark[2]}
\maketitle
\begin{abstract}
In 1991, one of the authors showed the existence of quadratic transformations
between the Painlev\'e VI equations with local monodromy differences
$(1/2,a,b,\pm 1/2)$ and $(a,a,b,b)$. In the present paper we give concise
forms of these transformations. 
They are related to the 
quadratic transformations obtained by Manin and
Ramani-Grammaticos-Tamizhmani via Okamoto transformations.
To avoid cumbersome expressions with differentiation,
we use contiguous relations instead of the Okamoto transformations.
The 1991 transformation is particularly important as it can be realized as
a quadratic-pull back transformation of isomonodromic Fuchsian equations.
The new formulas are illustrated by derivation of
explicit expressions for several complicated algebraic Painlev\'e VI functions.
\vspace{24pt}\\
{\bf 2000 Mathematics Subject Classification}: 34M55, 33E17. 
\vspace{24pt}\\
{\bf Short title}: {Quadratic transformations of Painlev\'e VI}\\
{\bf Key words}:  The sixth Painlev\'e equation, quadratic (or folding)
transformation, algebraic function.
\end{abstract}
\newpage
\setcounter{page}2

\section{Introduction}

The sixth Painlev\'e equation is, canonically,
\begin{eqnarray}
 \label{eq:P6}
\frac{d^2y}{dt^2}&=&\frac 12\left(\frac 1y+\frac 1{y-1}+\frac 1{y-t}\right)
\left(\frac{dy}{dt}\right)^2-\left(\frac 1t+\frac 1{t-1}+\frac 1{y-t}\right)
\frac{dy}{dt}\nonumber\\
&+&\frac{y(y-1)(y-t)}{t^2(t-1)^2}\left(\alpha+\beta\frac t{y^2}+
\gamma\frac{t-1}{(y-1)^2}+\delta\frac{t(t-1)}{(y-t)^2}\right),
\end{eqnarray}
where $\alpha,\,\beta,\,\gamma,\,\delta\in\CC$ are parameters. As well-known
\cite{JM}, its solutions
define isomonodromic deformations (with respect to $t$) of the $2\times2$
matrix Fuchsian equation with 4 singular points ($\lambda=0,1,t$, and
$\infty$):
\begin{equation}
 \label{eq:JM}
\frac{d}{d\lambda}\Psi=\left(\frac{A_0}\lambda+
\frac{A_1}{\lambda-1}+\frac{A_t}{\lambda-t}\right)\Psi,\qquad
\frac{d}{d\lambda}A_k=0\quad\mbox{for } k\in\{0,1,t\}.
\end{equation}
The standard correspondence is due to Jimbo and Miwa \cite{JM}. We choose
the traceless normalization of (\ref{eq:JM}), so we assume that the
eigenvalues of $A_0$, $A_1$, $A_t$ are, respectively, $\pm\theta_0/2$,
$\pm\theta_1/2$, $\pm\theta_t/2$, and that the matrix
$A_\infty:=-A_1-A_2-A_3$ is diagonal with the diagonal entries
$\pm\theta_\infty/2$.  Then the corresponding Painlev\'e equation has the
parameters
\begin{equation}
 \label{eq:para}
\alpha=\frac{(\theta_\infty-1)^2}2,\quad
\beta=-\frac{{\theta}_0^2}2,\quad\gamma=\frac{{\theta}_1^2}2,
\quad\delta=\frac{1-{\theta}_t^2}2.
\end{equation}
We refer to the numbers $\theta_0$, $\theta_1$, $\theta_t$ and
$\theta_\infty$ as {\em local monodromy differences}. They are invariants of
the isomonodromic deformation.

For any numbers $\nu_0,\nu_1,\nu_t,\nu_\infty$, let us denote by
$P_{VI}(\nu_0,\nu_1,\nu_t,\nu_\infty;t)$ the Painlev\'e VI equation for the
local monodromy differences $\theta_i=\nu_i$ for $i\in\{0,1,t,\infty\}$, via
(\ref{eq:para}). Note that changing the sign of $\nu_0,\nu_1,\nu_t$ or
$1-\nu_\infty$ does not change the Painlev\'e equation. 
Fractional-linear transformations for the Painlev\'e VI equation permute
the 4 singular points of (\ref{eq:JM}) and the numbers $\nu_0,\nu_1,\nu_t,1-\nu_\infty$.

The subject of this paper is quadratic transformations for the sixth
Painlev\'e equation. Their existence  was discovered in \cite{K5a},
\cite{K5}. In particular \cite{K5}, quadratic transformations were found
between isomonodromic Fuchsian equations (\ref{eq:JM}) with the local
monodromy differences $(\theta_0,\theta_1,\theta_t,\theta_\infty)$ related
as follows:
\begin{equation} \label{qua:kitaev}
\textstyle \left(a,a,b,b\right) \mapsto \left(\frac12,a,b,\frac12\right).
\end{equation}
These transformations act on the fundamental solution of (\ref{eq:JM}) as
$\Psi(\lambda)\mapsto S(\lambda)\Psi(R(\lambda))$, where $R(\lambda)$ is a
scalar quadratic function and $S(\lambda)$ is a matrix-valued rational
function. But corresponding transformation between Painlev\'e VI solutions
is implied as cumbersome compositions of lengthy formulas. This paper
presents compact expressions for quadratic transformation (\ref{qua:kitaev})
of Painlev\'e VI functions, up to fractional-linear transformations.

Simpler quadratic transformations for Painlev\'e VI equations are obtained
in \cite{M} and \cite{RGT}. Manin found that Landen's transformation for the
elliptic form of the Painlev\'e VI equation changes the local monodromy
differences as follows:
\begin{equation} \label{qua:manin0} \textstyle
(B,0,0,C)\mapsto\left(\frac{B}2,\frac{C-1}{2},\frac{B}2,\frac{C+1}2\right).
\end{equation}
The transformation in \cite{RGT} is the same, up to fractional-linear
transformations, as the rational form of Manin's transformation.
The local monodromy differences are changed as follows:
\begin{equation} \label{qua:manin} \textstyle
(0,A,B,1)\mapsto\left(\frac{A}2,\frac{B}{2},\frac{B}2,\frac{A}2+1\right).
\end{equation}
Explicit formulation of this transformation is very simple;  we present it in Lemma
\ref{lem:manin} below. To relate the transformations in (\ref{qua:kitaev}) and
(\ref{qua:manin}), set $C=A+1$.

The Painlev\'e equations involved in (\ref{qua:kitaev})--(\ref{qua:manin}), and
corresponding Fuchsian equations (\ref{eq:JM}) have the following properties:
\begin{enumerate}
\item[\it (i)] As mentioned above, transformation  (\ref{qua:kitaev}) is realized in \cite{K5}
as a quadratic pull-back transformation of corresponding Fuchsian equations (\ref{eq:JM}).
As a consequence, this transformation preserves finiteness (or infiniteness) of
the monodomy group of corresponding Fuchsian equations.
\item[\it (ii)] The Fuchsian equations corresponding to solutions of $P_{VI}(0,A,B,1;t)$ have
logarithmic singularities at $\lambda=0$ and $\lambda=\infty$,
except for the degenerate solution $y(t)\equiv0$.
(Recall that the Painlev\'e solution $y(t)$ is a rational multiple of the lower-left element
of $A_0$ by the Jimbo-Miwa correspondence. The point $\lambda=0$
is non-singular only if $A_0=0$.
The point $\lambda=\infty$ is logarithmic by the fractional-linear symmetry
$\lambda\mapsto1/\lambda$; the Painlev\'e solution is transformed as $y\mapsto1/y$.)
\item[\it (iii)] The Fuchsian solutions corresponding to solutions of
$P_{VI}\left(\frac{A}2,\frac{B}{2},\frac{B}2,\frac{A}2+1;t\right)$
do not have logarithmic singularities in general. Note incidentally that
the Painlev\'e equation has a simple solution $y(t)=\sqrt{t}$ for arbitrary $A, B$.
\end{enumerate}
It is apparent that quadratic transformations (\ref{qua:kitaev}) and (\ref{qua:manin})
have different character on the level of Fuchsian equations (\ref{eq:JM}). In particular,
contrary to {\it (i)} above,
transformation (\ref{qua:manin}) cannot be realized as a quadratic transformation of
the Fuchsian equations, because only one side necessarily has logarithmic points.
It is even possible that the monodromy group of equation (\ref{eq:JM})
on the $\left(\frac{A}2,\frac{B}{2},\frac{B}2,\frac{A}2+1\right)$ side of (\ref{qua:manin}) is finite,
whereas of course the monodromy group  cannot be finite on the $(0,A,B,1)$ side.
\begin{example} \rm \label{hitchinex}
Let us set $A=1$, $B=1$; then we have Hitchin's equation
$P_{VI}\left(\frac12,\frac12,\frac12,\frac12;t\right)$
on the $\left(\frac{A}2,\frac{B}{2},\frac{B}2,\frac{A}2+1\right)$ side. 
It is shown in \cite{Hit} (and in \cite{Hi} as well) that 
this equation has infinitely many algebraic solutions,
and they correspond to Fuchsian equations (\ref{eq:JM}) with finite dihedral
monodromy groups. Here is a parametrization of one such solution:
\begin{equation}
y(t)=-s, \qquad t=\frac{s^3(s+2)}{2s+1}.
\end{equation}
This is a fractional-linear version of the solution given at the end of Section 9 in \cite{Hit}.
The monodromy group of the corresponding Fuchsian equation 
is the dihedral group  with 6 elements.
Using formulas for 
quadratic transformation (\ref{qua:manin}) one computes that the corresponding solution
of $P_{VI}\left(0,A,B,1;\,\widehat{t}=(\sqrt{t}-1)^2/(\sqrt{t}+1)^2\right)$ can be parametrized as
\begin{equation}
\widehat{y}\left(\,\widehat{t}\,\right)=\frac{u(2u+1)}{u+2}, \qquad \widehat{t}=\frac{u^3(u+2)}{2u+1}.
\end{equation}
The parameters $u$, $s$ are related by the equation $(u+1)^2(s+1)^2+2us=0$.
Since $\widehat{y}\left(\widehat{t}\right)\not\equiv0$, the monodromy group of the
corresponding Fuchsian system is not finite.
(We continue to consider Hitchin's solutions in Section \ref{sec:prelim} and
Example \ref{ex:hitchin2}.)
\end{example}

Quadratic transformations (\ref{qua:kitaev}) and (\ref{qua:manin}) are
related by Okamoto transformations; this was noticed in \cite{CC},
\cite{GF}. An Okamoto transformation acts on the local monodromy differences
of Painlev\'e VI equations as follows:
\begin{equation} \label{eq:oka} 
\left( \theta_0,\,\theta_1,\,\theta_t,\,\theta_\infty \right) \mapsto
\left(\theta_0-\Theta,\,\theta_1-\Theta,\,\theta_t-\Theta,\,\theta_\infty-\Theta\right),
\end{equation}
where $\Theta=(\theta_0+\theta_1+\theta_t+\theta_\infty)/2$. In particular,
Okamoto transformations directly relate
\begin{equation} \label{eq:oka1} \textstyle
\left(a,a,-b,b\right)\;\mapsto\;\left(0,0,a+b,b-a\right),
\end{equation}
and
\begin{equation} \label{eq:oka2} \textstyle
\left(-\frac12,a,-b,\frac12\right)\;\mapsto\;
\left(\frac{b-a-1}2,\frac{a+b}2,-\frac{a+b}2,\frac{b-a+1}2\right).
\end{equation}
Recall again that changing the sign of the local monodromy differences
$\theta_0$, $\theta_1$, $\theta_t$ and $1-\theta_\infty$ does not change the
Painlev\'e VI equation, thus several Okamoto transformations can be performed
on the same Painlev\'e VI equation; see Lemma \ref{lem:okam} below.
The right-hand  sides in (\ref{eq:oka1})--(\ref{eq:oka2}) are related
by (\ref{qua:manin0}) and fractional-linear transformations.

Although quadratic transformation (\ref{qua:kitaev}) is related to the
simpler transformation (\ref{qua:manin}) via a couple of Okamoto transformations,
it is useful to have a direct formula for (\ref{qua:kitaev}). Application of the quadratic
transformations to algebraic Painlev\'e VI functions in Sections \ref{sec:algebraic}
and \ref{sec:morealg} illustrates this handily. The different nature of transformations
(\ref{qua:kitaev}) and (\ref{qua:manin}) on the level of Fuchsian equations (\ref{eq:JM})
becomes clear as well. In particular, Okamoto transformations can change
the monodromy group of the corresponding Fuchsian equation (\ref{eq:JM}).
In our main examples, the monodromy group of Fuchsian equations on the $(a,a,b,b)$
and $\left(\frac12,\frac12,a,b\right)$ levels is finite --- specifically, the icosahedral
group. But the monodromy group 
on the $(0,A,B,1)$ level is certainly not finite.

A practical advantage of our formulas is that we avoid cumbersome composition
of algebraic and differential transformations, as it is the case with direct composition
of (\ref{qua:manin}) with Okamoto transformations.  It appears that there is no direct
algebraic relation between one solution of $P_{VI}(a,a,b,b;t)$ and one solution of
$P_{VI}\left(\frac12,a,b,\frac12\right)$. Therefore we present our new formulas as
extended {\em\bf contiguous relations}. Like usual contiguous relations for
Painlev\'e VI or Gauss hypergeometric functions, the new formulas relate three functions:
one function from the $\left(\frac12,a,b,\frac12\right)$ side and two functions from
the $(a,a,b,b)$ side, or vice versa. Contiguous relations are consequences of
Okamoto 
transformations, but they do not involve differentiation.

Our original motivation for this work was to provide new examples of algebraic
Painlev\'e VI functions. In particular, we were interested in the algebraic functions
corresponding to Fuchsian systems (\ref{eq:JM}) with the icosahedral
monodromy group. These functions were classified by Boalch in \cite{Bo2}; there
are 52 classes, reminiscent to the 15 Schwartz classes of algebraic hypergeometric functions.
In the sixth electronic version of \cite{Bo2}, ten Boalch classes were not exemplified
yet\footnote{Currently, all 52 Boalch
icosahedral classes of Painlev\'e VI functions are exemplified; see the final version of
\cite{Bo2}. Of the mentioned 10 examples, we independently computed a type
42 example as well \cite{KV1}.
}. Eight of the missing examples can be obtained from earlier known by quadratic transformations.
The 8 examples were computed at the same time by Boalch \cite{Bo4} and us.
Compared to \cite{Bo4}, we derive the 8 examples most conveniently by
employing our new direct formulas for quadratic transformation (\ref{qua:kitaev}).
We are able to present explicit expressions for Painlev\'e VI solutions
on the most motivating $\left(\frac12,\frac12,a,b\right)$ level\footnote{The formulas in
\cite{Bo4} give algebraic Painlev\'e VI functions on the
$\left(\frac{A}2,\frac{B}{2},\frac{B}2,\frac{A}2+1\right)$ level, but
updated electronic version of \cite{Bo4} is supplemented by {\sf Maple} code
with examples on the $\left(\frac12,\frac12,a,b\right)$ level as well.} .

The authors are thankful to Yousuke Ohyama for the invitation to 14th International
Summer School on Functional Equations in Okayama (Japan), August 10--13 (2005),
where we presented our results.

\section{The old and new results}
\label{sec:oldnew}

Here we review basic results on quadratic transformation (\ref{qua:manin}),
the Okamoto transformations and contiguous relations. Then we present our main
results: compact formulas for Kitaev's original transformation (\ref{qua:kitaev}).

First we note that recently, in \cite{TOS}, a general notion of {\em folding transformations}
for Painlev\'e equations is introduced. These transformations correspond to
fixed points on the space of local monodromy differences 
of B\"acklund transformations induced by Cremona isometries.
The quadratic transformations for the Painlev\'e VI
equation are instances of folding transformations.
(For readers familiar with \cite{TOS}, in Appendix Section \ref{sec:other}
we briefly explain our approach and results in the notation of \cite{TOS}.)

Here is the explicit formulation \cite{RGT} of quadratic transformation (\ref{qua:manin}).
\begin{lemma} \label{lem:manin}
Suppose that $y_1$ is a solution of $P_{VI}(0,A,B,1\,;\,t_1)$. Let us denote
\begin{equation} \label{def:taueta}
\tau=\sqrt{t_1},\qquad \eta=\sqrt{y_1},\qquad
T_1=\frac{(\tau+1)^2}{(\tau-1)^2}.
\end{equation}
Then the function
\begin{equation}
Y_1(T_1)=\frac{(\tau+1)(\eta+1)}{(\tau-1)(\eta-1)}
\end{equation}
is a solution of
$P_{VI}\left(\frac{A}2,\frac{B}{2},\frac{B}2,\frac{A}2+1;\,T_1\right)$.
\end{lemma}
\begin{proof} \rm The claim can be checked by direct computations.
\end{proof}

To have a convenient notation for Okamoto transformations, we introduce the
following operator on functions. For any
$\nu_0,\nu_1,\nu_t,\nu_\infty\in\CC$, let
\begin{equation} \label{def:oka}
K_{[\nu_0,\nu_1,\nu_t,\nu_\infty;\,t]}\,y(t):=
y(t)+\frac{\nu_0+\nu_1+\nu_t+\nu_\infty}{Z(t)},
\end{equation}
where
\begin{equation} \label{def:okaz}
Z(t)=\frac{(t-1)\,\frac{dy(t)}{dt}-\nu_0}{y(t)}
-\frac{t\,\frac{dy(t)}{dt}+\nu_1}{y(t)-1}
+\frac{\frac{dy(t)}{dt}-1-\nu_t}{y(t)-t}.
\end{equation}
Okamoto's result in \cite{O} can be formulated as follows.
\begin{lemma} \label{lem:okam}
Suppose that $y(t)$ is a solution of
$P_{VI}(\theta_0,\theta_1,\theta_t,\theta_\infty;\,t)$, and that
$$\nu_0\in\{\theta_0,-\theta_0\},\quad \nu_1\in\{\theta_1,-\theta_1\},\quad
\nu_t\in\{\theta_t,-\theta_t\},\quad
\nu_\infty\in\{\theta_\infty,2-\theta_\infty\}.$$ Let
$$\Theta=\frac{\nu_0+\nu_1+\nu_t+\nu_\infty}2.$$ Then the function
$K_{[\nu_0,\nu_1,\nu_t,\nu_\infty;\,t]}\,y(t)$ is a solution of
$P_{VI}(\nu_0-\Theta,\nu_1-\Theta,\nu_t-\Theta,\nu_{\infty}-\Theta;t)$.
\end{lemma}
\begin{proof} \rm See \cite{O}.
The claim can be checked by direct computations.
\end{proof}

From a single Painlev\'e VI equation one can get up to 16 different
Painlev\'e VI equations by Okamoto transformations.
As shown in \cite{O}, subsequent combinations of Okamoto transformations
acting on the space of tuples $(\theta_0,\theta_1,\theta_t,\theta_\infty)$
form a group isomorphic to the affine Weyl group of type $D_4$.
If  all fractional-linear transformations are taken into account as well,
the symmetry group is the affine Weyl group of type $F_4$.
Commutativity of Okamoto (and fractional-linear) transformations
within the $D_4$ and $F_4$ lattices is recaped in Lemmas \ref{lem:frlin}
and \ref{lem:okamcom} below. In our paper, transformations of
the $\theta$-tuples leaving the Painlev\'e VI equation invariant
are considered as trivial. 
But the symmetry group is still evidently related to $D_4$.
In particular, combinations of Okamoto transformations relate
pairs of Painlev\'e VI functions whose respective local monodromy
differences are shifted by (either all odd or all even) integers.
For a precise statement, see Lemma \ref{lem:schlesinger} in the Appendix Section.

Generally, one can take any three Painlev\'e VI functions related to each
other by chains of Okamoto transformations,
and write down a nonlinear relation between them without the derivatives of
those functions \cite{JM,FGR,CM}. (For example,
one may take two Okamoto transformations of one Painlev\'e VI function and
eliminate the derivative of that function.) We refer to these relations as
{\em contiguous relations}; they are analogous to contiguous relations for Gauss
hypergeometric functions. Contiguous relations are usually more compact than
differential expressions for subsequent Okamoto transformations.

We prove our main formulas by using contiguous relations between functions
in the Okamoto orbits of $(a,a,b,b)$ and of $(1/2,1/2,a,b)$.   Instead of
composing quadratic transformation (\ref{qua:manin}) with Okamoto
transformations (\ref{eq:oka1})--(\ref{eq:oka2}), we rather compose it with
contiguous relations. Effectively, we merge
the two contiguous orbits into one via quadratic relation (\ref{qua:manin}).
In this sense, our formulas are extended contiguous relations,
as mentioned in Introduction.

Here below we present our main results. Theorem \ref{the:main} presents
the fractional-linear version of (\ref{qua:kitaev}) where the points
$\lambda=0$ and $\lambda=\infty$ are symmetric.
With this symmetry, the proofs are least cumbersome, and the formulas
are most elegant. Theorem \ref{the:main} presents quadratic transformation
(\ref{qua:kitaev}) directly. There the points $\lambda=0$ and $\lambda=1$ are
symmetric. Formulas with this symmetry are most convenient for applications
to algebraic Painlev\'e VI functions. To see the connection between both theorems,
note that $P_{VI}(a,a,b,b\,;t)$ and $P_{VI}(a,b-1,b,a+1\,;t/(t-1))$ are related by
a fractional-linear transformation.
\begin{theorem} \label{the:main}
Suppose that $y_0$ is a solution of $P_{VI}(a,b-1,b,a+1\,;t_1)$. Let us
denote
\begin{equation} \label{eq:okatr1}
y_1=K_{[-a,\,b-1,-b,\,1-a;\,t_1]}\,y_0,
\end{equation}
and
\begin{equation} \label{def:taueta2}
\tau=\sqrt{t_1},\qquad \eta=\sqrt{y_1},\qquad
T_1=\frac{(\tau+1)^2}{(\tau-1)^2}.
\end{equation}
Then the function
\begin{eqnarray} \label{eq2:sol1212}
Y_0(T_1)=\frac{a(\tau+1)(\eta+1)(y_0+\eta)}{(\tau-1)\left(a\eta(y_0-1)-(b-1)\,(y_0-\eta^2)\right)}
\end{eqnarray}
is a solution of $P_{VI}\left(a,\frac{1}2,\frac{1}2,b;T_1\right)$.
\end{theorem}

\begin{theorem} \label{lem:main}
Suppose that $g_0$ is a solution of $P_{VI}(a,a,b,b\,;t_2)$. Let us denote
\begin{equation} \label{eq2:okatr1}
g_1=K_{[-a,-a,-b,\,b\,;\,t_2]}\,g_0
\end{equation}
and
\begin{equation} \label{eq2:maint}
T_2=\frac12+\frac{t_2-\frac12}{2\,\sqrt{t_2^2-t_2}}.
\end{equation}
Then the function
\begin{eqnarray} \label{eq2:mainsol}
G_0(T_2)=\frac12+\frac{t_2-g_1+\sqrt{g_1^2\!-\!g_1}}{2\,\sqrt{t_2^2-t_2}}
+\frac{(a-b+1)(g_0-g_1)\left(g_1-\frac12-\sqrt{g_1^2\!-\!g_1}\right)}
{2\left(a\sqrt{g_1^2-g_1}-(b-1)(g_0-g_1)\right)\sqrt{t_2^2-t_2}}
\end{eqnarray}
is a solution of $P_{VI}\left(\frac12,\frac12,a,b\,;\,T_2\right)$.
\end{theorem}
The branches of the functions $\sqrt{t_1}$, $\sqrt{y_1}$, $\sqrt{t_2^2-t_2}$, $\sqrt{g_1^2-g_1}$
in these two theorems can be chosen arbitrary. We prove them in Section \ref{sec:proof}.
In Section \ref{sec:varia} we present some variations of our formulas. Also, in Appendix Section
\ref{sec:other} we give the main results in an alternative notation.

The authors have prepared a {\sf Maple} worksheet illustrating the formulas of this article.
The worksheet accessed by contacting the first author, or consulting his current webpage.

\section{Preliminaries}
\label{sec:prelim}

Here we present some simple results and observations, which are useful in
our arguments. For the sake of complete picture, we also refer to Table \ref{fig:frli}
in Appendix Section \ref{sec:other} of all fractional-linear transformations for the
Painlev\'e VI equation and its solutions.

It is useful to note that Okamoto transformations commute with the
fractional-linear transformations, and that when they act on Painlev\'e
functions they commute according to the $D_4$ (or $F_4$) lattice.
\begin{lemma} \label{lem:frlin}
Suppose that $(\alpha,\beta,\gamma,\delta)$ is a permutation of
$(0,1,t,\infty)$, and let $L:(y,t)\mapsto (Y,T)$ with
$T\in\{t,1-t,t/(t-1),1/t,1/(1-t),(t-1)/t\}$ denote the corresponding
fractional-linear transformation, as in Table $\ref{fig:frli}$. Then for any
numbers $\nu_0,\nu_1,\nu_t,\nu_\infty\in\mathbb C$ we have
\begin{equation}
L\,K_{[\nu_0,\,\nu_1,\,\nu_t,\,1+\nu_\infty;\,t]}=
K_{[\nu_\alpha,\,\nu_\beta,\,\nu_\gamma,\,1+\nu_\delta;\,T]}\,L.
\end{equation}
\end{lemma}
\begin{proof} It is enough to check the statement explicitly for a
generating set of the permutations. One can take, for example, the three
transpositions realized by the substitutions $\lambda\mapsto 1-\lambda$,
$\lambda\mapsto 1/\lambda$, $\lambda\mapsto t\lambda$ as in Table
\ref{fig:frli}.
\end{proof}
\begin{lemma} \label{lem:okamcom}
Suppose that $y(t)$ is a solution of
$P_{VI}(\nu_0,\nu_1,\nu_t,\nu_\infty;\,t)$, and let
$\Theta$ denote $(\nu_0+\nu_1+\nu_t+\nu_\infty)/2$. 
Then
\begin{equation}
K_{[\nu_0-\Theta,\,\nu_1-\Theta,\,\Theta-\nu_t,\,\nu_\infty-\Theta;\,t]}
K_{[\nu_0,\,\nu_1,\,\nu_t,\,\nu_\infty;\,t]}\,y(t)=
K_{[\nu_0,\,\nu_1,\,-\nu_t,\,\nu_\infty;\,t]}\,y(t)
\end{equation}
Besides,
\begin{equation} \label{eq1:oinvert}
K_{[\nu_0-\Theta,\,\nu_1-\Theta,\,\nu_t-\Theta,\,\nu_\infty-\Theta;\,t]}
K_{[\nu_0,\,\nu_1,\,\nu_t,\,\nu_\infty;\,t]}\,y(t)=y(t).
\end{equation}
\end{lemma}
\begin{proof} \rm The statements can be checked by direct computations. (The latter claim
is a convenient equivalent of the Painlev\'e VI equation.)
\end{proof}

Fractional-linear versions of quadratic transformation (\ref{qua:manin}) are
concisely presented in Table \ref{fig:quadr}. Extending the statement of Lemma
\ref{lem:manin},  this table 
can be used to compute any fractional-linear version of quadratic transformation
(\ref{qua:manin}) as follows. One may start with a
Painlev\'e VI solution $y(t)$ represented by one of the first six rows of
the table, compute $\tau$ and $\eta$ from the given expressions for $y$ and
$t$, and then pick up one of the bottom three rows, read off local monodromy
differences and an expression in terms of $\tau$ and $\eta$ of other
Painlev\'e VI solution. Or one may go the other direction: start with one of
the bottom three rows and get a transformed function for one of the six top
rows. Table \ref{fig:quadr} may be extended to include entries with $A$
interchanged with $B$; in the extra entries $\eta$ should be replaced by $\tau/\eta$.
\begin{table}\begin{center}
\begin{tabular}{|c|c|c|c|c|c|} \hline
$(\theta_0,\theta_1,\theta_t,\theta_{\infty})$
 & $y$ & $t$ \\ \hline
$(0,A,B,1)$ & $\eta^2$ & $\tau^2$ \\
$(A,0,B,1)$ & $1-\eta^2$ & $1-\tau^2$ \\
$(A,B,0,1)$ & $(\eta^2-1)/(\tau^2-1)$ & $1/(1-\tau^2)$ \\
$(B,0,0,A+1)$ & $(\eta^2-\tau^2)/(\eta^2-1)$ & $\tau^2$ \\
$(0,B,0,A+1)$ & $(\tau^2-1)/(\eta^2-1)$ & $1-\tau^2$ \\
$(0,0,B,A+1)$ & $1/(1-\eta^2)$ & $1/(1-\tau^2)$ \vspace{-4pt} \\
\dotfill & \dotfill & \dotfill \\
$(\frac{A}2,\frac{B}2,\frac{B}2,\frac{A}2+1)$ & $(\tau+1)(\eta+1)/(\tau-1)(\eta-1)$ & $(\tau+1)^2/(\tau-1)^2$ \\
$(\frac{A}2,\frac{A}2,\frac{B}2,\frac{B}2+1)$ & $(\tau+1)(\eta+1)/2(\eta+\tau)$ & $(\tau+1)^2/4\tau$ \\
$(\frac{B}2,\frac{A}2,\frac{B}2,\frac{A}2+1)$ & $2(\eta+\tau)/(\tau+1)(\eta+1)$ & $4\tau/(\tau+1)^2$ \\
 \hline
\end{tabular}\end{center}
\caption{Simple quadratic transformations} \stepcounter{figure}
\label{fig:quadr}
\end{table}

Within setting of Lemma \ref{lem:manin}, let us compare the function fields
$\CC(t_1,y_1)$ and $\CC(T_1,Y_1)$. At first glance, the quadratic
transformation requires to adjoin $\tau$ and $\eta$ to $\CC(t_1,y_1)$ in
order to get $\CC(T_1,Y_1)$. However, the automorphism $\tau\mapsto 1/\tau$,
$\eta\mapsto1/\eta$ fixes the field $\CC(T_1,Y_1)$. We have the following
diagram, where all immediate field extensions have degree 2 in general:
\begin{equation} \label{eq:fields}
\begin{array}{ccccl} \CC(t_1,y_1) &\subset & \CC\left(\tau,y_1\right) &
\subset & \CC\left(\tau,\eta\right) \\
&& \cup && \quad\cup\\
&& \CC\left(T_1,y_1\right) & \subset & \CC\left(T_1,Y_1\right). 
\end{array}
\end{equation}
As we see, $\CC(T_1,Y_1)$ is generally an index 2 subfield of a degree 4 extension of
$\CC(t_1,y_1)$. In particular, if $y_1$ is an algebraic function, then the
algebraic degree of the extension $\CC(T_1,Y_1)\supset\CC(T_1)$ is usually
twice the degree of $\CC(t_1,y_1)\supset\CC(t_1)$.
However, Example \ref{hitchinex} gives an explicit situation when
algebraic solutions on both sides have the same degree, four.
In this case, $\CC\left(\tau,y_1\right)=\CC\left(\tau,\eta\right)$.

Algebraic geometrically, the quadratic transformation is a 4-to-2
correspondence in general: the projection onto the $(\tau,\eta)$-plane relates 4
analytic branches on the $(\frac{A}2,\frac{B}2,$$\frac{B}2,\frac{A}2\!+\!1)$
side with 2 branches on the $(0,A,B,1)$ side. Permutation of the former
branches is realized by $\tau\to-\tau$ and/or $\eta\to-\eta$, and also by
fractional-linear permutations of the two $\frac{A}2$ points and (or) the
two $\frac{B}2$ points. Permutation of the latter branches is realized by
simultaneous $\tau\mapsto 1/\tau$, $\eta\mapsto1/\eta$, and also by
fractional-linear permutation of the $0/1$-points.

\section{Proof of the main results}
\label{sec:proof}

An important intermediate question for us is when the property of being related by
Okamoto transformations is preserved by the quadratic transformation.
Notation throughout the paper is entirely
consistent if we identify $A=a+b-1$, $B=b-a$.

\begin{lemma} \label{lem:okkam}
Suppose that $y_1$ is a solution of $P_{VI}(0,A,B,1;t_1)$, and that $y_2$ is
a solution of $P_{VI}(0,B-1,A+1,1;t_1)$. Suppose that
\begin{equation} \label{def:g0g0}
K_{[0,\,A,-B,\,1;\,t_0]}\,y_1=K_{[0,\,B-1,\,-A-1,\,1;\,t_0]}\,y_2.
\end{equation}
Let $y_0$ denote the evaluation of any side of this equality. Let $\tau$,
$\eta$, $T_1$ be defined as in $(\ref{def:taueta2})$.
Then the functions
\begin{equation} \label{def:y0y2}
Y_1(T_1)=\frac{(\tau+1)(\eta+1)}{(\tau-1)(\eta-1)}, \qquad
Y_2(T_1)=\frac{(\tau+1)(y_0+\eta)}{(\tau-1)(y_0-\eta)}
\end{equation}
are solutions of, respectively,
\begin{equation} \label{eq2:painls} \textstyle
P_{VI}\left(\frac{A}2,\frac{B}2,\frac{B}2,\frac{A}{2}+1;\,T_1\right)\qquad
\mbox{and}\qquad
P_{VI}\left(\frac{B-1}2,\frac{A+1}{2},\frac{A+1}2,\frac{B+1}2;\,T_1\right),
\end{equation}
and we have
\begin{equation} \label{eq2:qoka}
Y_2=K_{\left[\frac{A}2,-\frac{B}2,-\frac{B}2,\frac{A}2+1;\,T_1\right]}Y_1.
\end{equation}
\end{lemma}
\begin{proof} The contiguous relation between $y_1$, $y_2$, $y_0$ can be
derived by expressing $y_1$, $y_2$ as Okamoto transformations of $y_0$, and
eliminating the derivative of $y_0$ from the two identities. The result is
very simple:
\begin{equation} \label{eq1:back}
y_1y_2=y_0^2.
\end{equation}
By $\eta$ we actually denote a branch of $\sqrt{y_1}$.  We have two choices
$\sqrt{y_2}=\pm y_0/\eta$; we choose $\sqrt{y_2}=y_0/\eta$. Now we apply
Lemma \ref{lem:manin} to $y_1$ and $y_2$,
and conclude that the functions $Y_1$ and $Y_2$ satisfy respective
Painlev\'e VI equations in (\ref{eq2:painls}).
Formula (\ref{eq2:qoka}) depends on the right choice of $\sqrt{y_2}$ we did.
To show that formula, we express $dy_1/dt$ in terms of $y_1,y_0$ by using the
definition of $y_0$ by the left-hand side of (\ref{def:g0g0}). Then we
easily express $d\eta/d\tau$ in terms of $\eta,y_0$:
\begin{equation} \label{eq2:dedt}
\frac{d\eta}{d\tau}=\frac{\eta\left(A(\eta^2-\tau^2)(y_0-1)
+(1-B)(\eta^2-1)(y_0\!-\!\tau^2)\right)}{\tau(\tau^2-1)(y_0-\eta^2)}.
\end{equation}
 Expression (\ref{eq2:qoka}) can be rewritten in
terms of $d\eta/d\tau$, $\eta,y_0$ by using (\ref{def:y0y2}). After
substituting (\ref{eq2:dedt}) we check the identity.
\end{proof}

In the above Lemma, note that $y_0$ satisfies
$P_{VI}\left(\frac{B-A-1}2,\frac{A+B+1}2,\frac{A+B+1}2,\frac{B-A+1}2;t_1\right)$.
If we would replace $y_0\mapsto-y_0$ in (\ref{def:y0y2}), the function $Y_2$
would still be a solution of
$P_{VI}\!\left(\frac{B-1}2,\frac{A+1}{2},\frac{A+1}2,\frac{B+1}2;\,T_1\right)$,
but identity(\ref{eq2:qoka}) would not hold.

Now we are ready to prove the main results.
As mentioned above, our strategy to relate the Painlev\'e VI equations with
local monodromy differences $(a,a,b,b)$ and $(\frac12,a,b,\frac12)$
is to combine the quadratic relation of  Lemma \ref{lem:manin} with contiguous relations
of Painlev\'e VI equations on both sides. A contiguous relation on the $(a,a,b,b)$ side
is indirectly employed by a reference to Lemma \ref{lem:okkam}; that contiguous relation
is (\ref{eq1:back}).


\vspace{6pt}
{\bf\noindent Proof of Theorem \ref{the:main}.}
We assumed that $y_0$ is a solution of $P_{VI}(a,b-1,b,a+1\,;t_1)$. We
defined $y_1\!=K_{[-a,\,b-1,-b,\,1-a;\,t_1]}\,y_0$, hence $y_1$ is a
solution of $P_{VI}(0,a+b-1,b-a,1;t_1)$. The variables $\tau$, $\eta$ and
$T_1$ are defined as in $(\ref{def:taueta2})$.

Let us denote $y_2=K_{[a,\,b-1,-b,\,a+1;\,t_1]}\,y_0$. It is a solution of
\mbox{$P_{VI}(0,b-a-1,a+b,1;t_1)$}. Let $Y_1$ and $Y_2$ be defined as in formula
(\ref{def:y0y2}). Lemma \ref{lem:okkam} tells us that $Y_1$ and $Y_2$ are
quadratic transformations of $y_1$ and $y_2$, respectively, and
\begin{equation} \label{eq2:qokay}
Y_2=K_{\left[\frac{a+b-1}2,\,\frac{a-b}2,\,\frac{a-b}2,\,\frac{a+b+1}2;\,T_1\right]}\,Y_1.
\end{equation}
Let us consider
\begin{equation}
Y_0=K_{\left[\frac{1-a-b}2,\,\frac{a-b}2,\,\frac{a-b}2,\,\frac{a+b+1}2;\,T_1\right]}\,Y_1.
\end{equation}
This is a solution of $P_{VI}\left(a,\frac{1}2,\frac{1}2,b;T_1\right)$. The
contiguous relation between $Y_1,Y_2$ and $Y_0$ can be computed similarly as
(\ref{eq1:back}). The result is
\begin{eqnarray} \label{eq3:backlund}
Y_0\equal\frac{2aY_1Y_2}{(a-b+1)Y_1+(a+b-1)Y_2}.
\end{eqnarray}
Rewriting $Y_0$ in terms of $\tau$, $\eta$, $y_0$ gives (\ref{eq2:sol1212}).
\hfill$\Box$\bigskip

{\bf\noindent Proof of Theorem \ref{lem:main}.}
This is a fractional-linear version of Theorem \ref{the:main}.
We identify:
\begin{equation} \label{eq:frlis}
t_2=\frac{t_1}{t_1-1}, \quad
g_0=\frac{y_0}{y_0-1}, \quad g_1=\frac{y_1}{y_1-1}, \quad
T_2=\frac{T_1}{T_1-1}, \quad
G_0=\frac{Y_0-T_1}{1-T_1}.
\end{equation}
Consequently, we choose
\begin{equation} \label{eq:tauet}
\tau=\frac{\sqrt{t_2^2-t_2}}{t_2-1},\qquad
\eta=-\frac{\sqrt{g_1^2-g_1}}{g_1-1}.
\end{equation}
Then we substitute this to (\ref{eq2:sol1212}).
\hfill$\Box$\bigskip

Theorem \ref{lem:main} can be proved directly, without reference to Theorem \ref{the:main},
but using the same proof scheme. 
In particular, one may consider $g_2=K_{[a,a,-b,b;\,t_2]}\;g_0$,
and the following quadratic transformations of $g_1$ and $g_2$,
in the common terms of Table \ref{fig:quadr} and (\ref{eq:frlis})--(\ref{eq:tauet}):
\begin{equation}
G_1=\frac{(\tau+1)(\eta-\tau)}{2\tau(\eta-1)},\qquad\qquad
G_2=\frac{(\tau+1)(y_0-\tau\eta)}{2\tau(y_0-\eta)}.
\end{equation}
These are solutions of
$P_{VI}\!\left(\frac{b-a}2, \frac{b-a}2, \frac{a+b-1}2, \frac{a+b+1}2;T_2\!\right)$ and
$P_{VI}\!\left(\frac{a+b}2, \frac{a+b}2, \frac{b-a-1}2, \frac{b-a+1}2;T_2\!\right)$,
respectively. We would have
\begin{equation} \label{pprf:oka1}
G_2 = K_{\left[ \frac{a-b}2, \frac{a-b}2, \frac{a+b-1}2, \frac{a+b+1}2;\,
T_2\right]}\,G_1\quad\mbox{and}\quad G_0 = K_{\left[ \frac{a-b}2,
\frac{a-b}2, \frac{1-a-b}2, \frac{a+b+1}2;\, T_2\right]}\,G_1.
\end{equation}
However, the contiguous relations between $g_0,g_1,g_2$ and $G_0,G_1,G_2$
are more messy than (\ref{eq1:back}) and (\ref{eq3:backlund}). The
intermediate expressions are:
\begin{eqnarray}
\label{prf:RQY} G_1(T_2)\equal \frac12+\frac{t_2-g_1+\sqrt{g_1^2-g_1}}{2\,\sqrt{t_2^2-t_2}}, \\
\label{prf:RQQ} G_2(T_2)\equal
\frac12+\frac{t_2-g_2+\sqrt{g_2^2-g_2}}{2\,\sqrt{t_2^2-t_2}},
\end{eqnarray}
where
$\sqrt{g_2^{2}-g_2}=-(g_0^2-g_0)\sqrt{g_1^2-g_1}\big/\left(g_0^2-2g_1g_0+g_1\right)$.
This relation between two square roots is a more complicated equivalent of
Lemma \ref{lem:okkam}.

The relation between the function fields $\CC(y_0,t_1)$ and $\CC(Y_0,T_1)$,
or between $\CC(g_0,t_2)$ and $\CC(G_0,T_2)$, is the same
as between the function fields $\CC(y_1,t_1)$ and $\CC(Y_1,T_1)$ in (\ref{eq:fields}),
because Okamoto and fractional-linear transformations do not change function fields.
In particular, formula (\ref{eq1:oinvert}) shows ``invertability" of Okamoto transformations.
The permutation $\tau\mapsto1/\tau,\eta\mapsto1/\eta$ of analytic branches on the
$(a,a,b,b)$ side should be supplemented with $y_0\mapsto1/y_0$. We have the
same branch permutations $\tau\mapsto -\tau$ and $\eta\mapsto -\eta$ on the
other side. However, it appears that the permutation $\eta\mapsto-\eta$ is not
realizable by fractional-linear transformations of Painlev\'e VI equations in the
$\left(\frac12,a,b,\frac12\right)$-orbit.

For the Painlev\'e VI equations $P_{VI}(a,a,b,b;t_2)$ and $P_{VI}(\frac12,\frac12,a,b;T_2)$,
the involution with $\tau\mapsto1/\tau$,  $\eta\mapsto1/\eta$ acts as the fractional-linear
transformation $\lambda\to1-\lambda$ (in terms of Table \ref{fig:frli} below).
In general, $\CC(G_0,T_2)$ must be an index 2 subfield of
$\CC(g_0,g_1,t_2,$ $\sqrt{t_2^2-t_2},\sqrt{g_1^2-g_1})$ fixed by this involution.
The involution acts as multiplication of $t_2-\frac12$, $g_1-\frac12$,
$g_0-g_1$, $\sqrt{t_2^2-t_2}$ and $\sqrt{g_1^2-g_1}$ by $-1$. Accordingly,
expression (\ref{eq2:mainsol}) for $G_0$ can be rewritten in the form
\begin{equation} \label{eq1:fmaint}
G_0(T_2)=\frac12+\widetilde{A}\,\frac{g_0-g_1}{\sqrt{t_2^2-t_2}}+
\widetilde{B}\,\frac{\sqrt{g_1^2-g_1}}{\sqrt{t_2^2-t_2}},
\end{equation}
with $\widetilde{A}$, $\widetilde{B}$ rational expressions in
$(t_2-\frac12)/(g_0-g_1)$, $(g_1-\frac12 )/(g_0-g_1)$ and $(g_1^2-g_1)/(g_0-g_1)^2$.
The form in (\ref{eq1:fmaint}) is often the most
convenient to represent ``quadratically" transformed algebraic Painlev\'e VI functions.

\section{Variations of new formulas}
\label{sec:varia}

We use the same functions $y_0,Y_0,Y_1,Y_2$ as in the previous section.
For shorthand convenience, let us denote by $\pi_2$ 
the fractional-linear transformation $y\mapsto t(y-1)/(y-t)$.
This is consistent with Table \ref{fig:frli} below.

Within the setting of Theorem \ref{the:main}, a solution of
$P_{VI}\left(a,\frac{1}2,\frac{1}2,b+1;T_1\right)$ is
\begin{eqnarray} \label{eq2:sol1212a}
\widetilde{Y}_0(T_1)=
\frac{a(\tau+1)(\eta+\tau)(y_0+\eta\tau)}{(1-\tau)\left(a\eta(y_0-\tau^2)+
b\tau(y_0-\eta^2)\right)}
\end{eqnarray}
To obtain this formula, one may apply the same Theorem \ref{the:main} to
the solution $t_1/y_0$ of $P_{VI}\left(a,-b,1-b,a+1;T_1\right)$, which is
a fractional linear transformation of $y_0$. In other words,
we just have to substitute $b\mapsto1-b$, $y_0\mapsto\tau^2/y_0$,
$\eta\mapsto\tau/\eta$ into (\ref{eq2:sol1212}). We checked explicit expressions
for all fractional-linear versions of quadratic transformation (\ref{qua:kitaev}),
and the formulas (\ref{eq2:sol1212}), (\ref{eq2:sol1212a}) appear to be most compact.
The expressions would be terribly cumbersome if expressed in terms of $y_0$ and
$dy_0/dt$ rather than in terms of $y_0$ and $\eta$; in particular, $\eta$ would have
to be replaced by the square root of  (\ref{def:oka}) according to (\ref{eq:okatr1}).

To get a solution of $P_{VI}\left(\frac{1}2,a,b,\frac{1}2;T_1\right)$, we
consider the fractional-linear transformation $\pi_2\widetilde{Y}_0$
of (\ref{eq2:sol1212a}).  
An explicit expression is
\begin{eqnarray}\label{eq2:sol1212b}
\frac{(\tau+1)\left(2a(\eta+\tau)(y_0+\eta)+(b-a)(\tau-1)(y_0-\eta^2)\right)}
{(\tau-1)\left(2a(\eta+\tau)(y_0-\eta)+(b-a)(\tau+1)(y_0-\eta^2)\right)}.
\end{eqnarray}
Along with (\ref{qua:kitaev}), transformation
$(a,a,b,b)\mapsto\left(\frac{1}2,a,b,-\frac{1}2\right)$ was 
considered in \cite{K5} as well. To compute a solution of
$P_{VI}\left(\frac{1}2,a,b,-\frac{1}2;T_1\right)$, we consider
\begin{eqnarray}
Y_3\equal K_{\left[\frac{1-a-b}2,\,\frac{a-b}2,\,\frac{b-a}2,\,\frac{a+b+1}2;\,T_1\right]}Y_1,\\
Y_4\equal
K_{\left[\frac{a+b}2,\,\frac{b-a+1}2,\,\frac{a-b+1}2,\,\frac{a+b}2;\,T_1\right]}Y_3.
\end{eqnarray}
These two functions can be computed by using contiguous relations, similarly as (\ref{eq1:back}).
The contiguous relations are:
\begin{eqnarray}
Y_3\equal\frac{\left(2aY_1Y_2+(b-a)T_1Y_1-(a+b)T_1Y_2\right)Y_1}
{Y_1^2+(2a-1)Y_1Y_2+(b-a-1)T_1Y_1-(a+b-1)T_1Y_2},\\
Y_4\equal\frac{\left(  (a-b)(T_1-1)(Y_3-Y_1)Y_3-Y_1Y_3^2
+(T_1+1)Y_3^2-2T_1Y_3+T_1Y_1\right) Y_3}
{(Y_3-Y_1)\!\left((a+b)(Y_3^2+T_1) -2(bT_1+a)Y_3+Y_3^2\right)
-Y_1Y_3^2+(T_1+1)Y_1Y_3-T_1Y_3}. \nonumber
\end{eqnarray}
The function $Y_4$ is a solution of
$P_{VI}\left(\frac{1}2,a,b,-\frac{1}2;T_1\right)$. However its expression in
terms of $\tau$, $\eta$, $y_0$ is lengthy. The corresponding expression for fractional-linear transformation $T_1(Y_4-1)/(T_1-1)Y_4$ is much shorter; this is a solution
of $P_{VI}\left(a,b,\frac{3}2,\frac{1}2;T_1/(T_1\!-\!1)\right)$:
\begin{gather*}
\frac{a(\tau+1)}{\tau}\,
\frac{2a\eta(\eta+1)(y_0-\tau^2)+(a\!+\!b\!-\!1)(\tau\!-\!1)(\eta\!-\!\tau)(y_0\!-\!\eta^2)
-2b\tau(\eta\!+\!1)(y_0\!-\!\eta^2)}
{2a(\eta-\tau)(y_0+\eta)+(a-b)(\tau+1)(y_0-\eta^2)} \times\hspace{-7pt}\\
\frac{2a(\eta+1)(y_0^2-\tau^2\eta^2)-(a+b)(\tau+1)(y_0+\tau\eta)(y_0-\eta^2)-
(\tau-1)(y_0-\tau\eta)(y_0-\eta^2)}
{4a^2(\eta^2\!-\!\tau^2)(y_0\!+\!\eta)^2+4a(a\!-\!b)(\eta\!+\!\tau^2)(y_0\!+\!\eta)(y_0\!-\!\eta^2)
-\left((a\!-\!b)^2\!-\!1\right)(\tau^2\!-\!1)(y_0\!-\!\eta^2)^2}.
\end{gather*}

To compute other similar solutions, it is convenient to know explicitly how
to shift the parameters $a$, $b$ by integers. These formulas may spare
cumbersome computations of contiguous relations.
\begin{lemma}
Suppose that we have an expression in terms of $\tau$, $\eta$, $y_0$ of a
solution of a Painlev\'e VI equation in the orbit of
$P_{VI}\!\left(a,a,b,b\,;t_0\right)$, $P_{VI}(0,a+b-1,b-a,1;t_0)$,
$P_{VI}\!\left(\frac{a+b-1}2,\frac{b-a}2,\frac{b-a}2,\frac{a+b+1}2;T_0\right)$
or $P_{VI}\!\left(\frac12,a,b,\frac12;T_0\right)$ under Okamoto and
fractional-linear transformations. To get a similar expression for a
solution of Painlev\'e VI equation with the parameters $a$, $b$ shifted or
interchanged, one may apply the following substitutions:
\begin{eqnarray*}
& b\mapsto b+1, & \eta\mapsto
\frac{\tau\eta\left(a(y_0-\tau^2)-b(y_0-\eta^2)\right)}
{a\eta^2(y_0-\tau^2)-b\tau^2(y_0-\eta^2)},\\
&&\!y_0\mapsto
\frac{\tau^2\!\left(a(y_0-\tau^2)-b(y_0-\eta^2)\right)
\left(a\eta^2(y_0-\tau^2)-by_0(y_0-\eta^2)\right)}
{\left(a\eta^2(y_0-\tau^2)-b\tau^2(y_0-\eta^2)\right)
\left(ay_0(y_0-\tau^2)-b\tau^2(y_0-\eta^2)\right)};\\
& b\mapsto b-1, & \eta\mapsto \frac{\tau\eta(a(y_0-1)+(b-1)(y_0-\eta^2))}
{a\eta^2(y_0-1)+(b-1)(y_0-\eta^2)},\\
&& \!y_0\mapsto
\frac{\tau^2\!\left(a(y_0\!-\!1)+(b\!-\!1)(y_0\!-\!\eta^2)\right)
\left(a\eta^2(y_0\!-\!1)+(b\!-\!1)y_0(y_0\!-\!\eta^2)\right)}
{\left(ay_0(y_0\!-\!1)+(b\!-\!1)(y_0\!-\!\eta^2)\right)
\left(a\eta^2(y_0\!-\!1)+(b\!-\!1)(y_0\!-\!\eta^2)\right)};
 \\ & a\leftrightarrow b, & \!y_0\mapsto
\frac{a\eta^2(y_0-\tau^2)-b\tau^2(y_0-\eta^2)}{a\,(y_0-\tau^2)-b\,(y_0-\eta^2)};\\
& b\mapsto 1-b, & \eta\mapsto\frac{\tau}{\eta},\quad
y_0\mapsto\frac{\tau^2}{y_0}.
\end{eqnarray*}
To shift $a\mapsto a\pm1$, one may consequently apply $a\leftrightarrow b$,
$b\mapsto b\pm 1$, $a\leftrightarrow b$.
\end{lemma}
\begin{proof} Let $L$ denote the fractional-linear transformation
$y\mapsto t_1/y$. It induces the transformation $b\mapsto 1-b$, as we already
noticed in the beginning of this section. 
Let
\begin{equation}
\widetilde{y}_0=K_{\left[-a,\,b-1,\,b,\,1-a;\,t_0\right]}\,y_0.
\end{equation}
The transformation $a\leftrightarrow b$ does not change $\eta$, and
transforms $y_0$ to $\widetilde{y}_0$. The expression for $\widetilde{y}_0$
can be obtained from the contiguous relation between $y_0,y_1,\widetilde{y}_0$,
similarly as (\ref{eq1:back}).

To compute $b\mapsto b+1$, we need to compute solutions of
$P_{VI}\left(0,a+b,b-a+1,1;t_0\right)$ and
$P_{VI}\left(a,b,b+1,a+1;t_0\right)$. We consider, respectively,
$$L\,K_{\left[b,\,a-1,\,-a,\,b+1;\,t_0\right]}\,\widetilde{y}_0\qquad\mbox{and}\qquad
L\,K_{\left[b,\,a-1,\,a,\,b+1;\,t_0\right]}\,\widetilde{y}_0.$$
Again, contiguous relations between $\widetilde{y}_0,y_0$ and each of these
functions give the formulas.

To compute $b\mapsto b-1$, we may compose $b\mapsto 1-b$, $b\mapsto b+1$,
$b\mapsto 1-b$. The statement about $a\mapsto a\pm 1$ is clear.
\end{proof}

As an example, one may check the relation $b\mapsto b+1$ between
(\ref{eq2:sol1212}) and (\ref{eq2:sol1212a}). The transformation $b\to 1-b$
effectively interchanges $A=a+b-1$ and $B=b-a$ as local monodromy
differences.

It is possible to derive relations between two functions in the
$\left(\frac12,a,b,\frac12\right)$-orbit and one functions from the
$(a,a,b,b)$ orbit. For example, one may eliminate $\eta$ from expressions
(\ref{eq2:sol1212}) and (\ref{eq2:sol1212a}). However, the relations are too
complicated for print. We only mention a relation between the functions $Y_1(T_1)$,
$Y_2(T_1)$ of the $\left(\frac{a+b-1}2,\,\frac{b-a}2,\,\frac{b-a}2,\,\frac{a+b+1}2\right)$ level
and the solution $y_0(t_1)$ of the $\left(a,a,b,b\right)$ level:
\begin{equation}
y_0(t_1)=\frac{(Y_1+\sqrt{T_1})(Y_2+\sqrt{T_1})}{(Y_1-\sqrt{T_1})(Y_2-\sqrt{T_1})},
\end{equation}
where we should identify $\sqrt{T_1}=(\tau+1)/(\tau-1)$.

\section{Algebraic Painlev\'e VI functions}
\label{sec:algebraic}

Here we demonstrate how our formulas for quadratic transformation
(\ref{qua:kitaev}) can be applied to compute new examples algebraic
Painlev\'e VI functions. We mainly concentrate on icosahedral algebraic
Painlev\'e VI functions, classified in \cite{Bo2}.
But first 
we develop Example \ref{hitchinex} further.

\begin{example} \label{ex:hitchin2} \rm
Quadratic transformation (\ref{qua:kitaev}) can be applied to Hitchin's solutions
of $P_{VI}\left(\frac12,\frac12,\frac12,\frac12;t\right)$ iteratively. Our new formulas
are directly suitable therefore. On the other hand, Hitchin's
equation can be transformed to Picard's \cite{Pic} equation
$P_{VI}(0,0,0,1;t)$ by an Okamoto transformation.
Lemma \ref{lem:manin} can be iteratively applied to Picard's equation.
Algebraic Picard's solutions are described in \cite{FuchsP} and \cite[Section 2]{Ma1}.
The algebraic solutions correspond to $n$-division points on a general elliptic curve,
so they are related to the modular curves $X_1(n)$. This relation is noticed for
Hitchin's case as well \cite{Hi}. Within this correspondence, the mentioned iterative
applications of (\ref{qua:kitaev}) or Lemma \ref{lem:manin} double $n$; this can
be easily seen from the elliptic form of (\ref{qua:manin0}) and the elliptic
form \cite[Theorem 1.4]{M} of $P_{VI}(0,0,0,1;t)$. More generally, there are
not only quadratic but also algebraic transformations of arbitrary degree $k$
for the Picard and Hitchin case of Painlev\'e VI equation.
For the Picard equation in the elliptic form,  these transformastions correspond
to the multiplication by $k$ on a general elliptic curve.
\end{example}

{\em Icosahedral} algebraic Painlev\'e VI functions are associated
to  Fuchsian systems (\ref{eq:JM}) with the icosahedral monodromy group.
As shown in \cite{Bo2}, up to Okamoto and fractional-linear transformations
there are 52 classes of icosahedral Painlev\'e VI functions. Quadratic
transformations relate some pairs of these classes, as recaped in Table
\ref{fig:icosaq}. In the first column, we identify the icosahedral classes
by the numbers in Boalch's classification \cite[Table 1]{Bo2}. In the next
two columns, we give representative tuples of local monodromy differences
for the transformed classes. In the last two columns, we give the
transformation of algebraic degree and genus of the algebraic Painlev\'e VI
functions.
\begin{table}\begin{center}
\begin{tabular}{|c|c|c|c|c|c|} \hline
\;Boalch's\; & \multicolumn{2}{|c|}{Local monodromy differences} &
\;Algebraic\; & \;Genus\; \\
\cline{2-3} cases & $(a,a,b,b)$ & $(1/2,1/2,a,b)$ & degree & \\
 \hline
$21\Rightarrow 28$ & $\left(1/5,1/5,2/5,2/5\right)$ & $(1/2,1/2,1/5,2/5)$ & $\;\,5\Rightarrow10$ & $0\Rightarrow0$ \\
$31\Rightarrow 44$ & \;$(1/5,1/5,1/5,1/5)$\; & \;$(1/2,1/2,1/5,1/5)$\; & $10\Rightarrow20$ & $0\Rightarrow1$ \\
$32\Rightarrow 45$ & $(2/5,2/5,2/5,2/5)$ & $(1/2,1/2,2/5,2/5)$ & $10\Rightarrow20$ & $0\Rightarrow1$ \\
$39\Rightarrow 47$ & $(1/3,1/3,1/5,1/5)$ & $(1/2,1/2,1/3,1/5)$ & $15\Rightarrow30$ & $1\Rightarrow2$ \\
$40\Rightarrow 48$ & $(1/3,1/3,2/5,2/5)$ & $(1/2,1/2,1/3,2/5)$ & $15\Rightarrow30$ & $1\Rightarrow2$ \\
$41\Rightarrow 49$ & $(1/3,1/3,1/3,1/3)$ & $(1/2,1/2,1/3,1/3)$ & $18\Rightarrow36$ & $1\Rightarrow3$ \\
$44\Rightarrow 50$ & $(1/2,1/2,1/5,1/5)$ & $(1/2,1/2,1/2,1/5)$ & $20\Rightarrow40$ & $1\Rightarrow3$ \\
$45\Rightarrow 51$ & $(1/2,1/2,2/5,2/5)$ & $(1/2,1/2,1/2,2/5)$ & $20\Rightarrow40$ & $1\Rightarrow3$ \\
$49\Rightarrow 52$ & $(1/2,1/2,1/3,1/3)$ & $(1/2,1/2,1/2,1/3)$ & $36\Rightarrow72$ & $3\Rightarrow7$ \\
 \hline
\end{tabular}\end{center}
\vspace{-8pt}\caption{Quadratic transformations of icosahedral Painlev\'e VI
functions} \label{fig:icosaq} \stepcounter{figure}
\end{table}

Quadratic transformations can be used to compute examples for the higher
degree classes of icosahedral Painlev\'e VI functions. Examples for the
classes 44--45, 47--52 are difficult to compute by other means. Our formulas
for quadratic transformation (\ref{qua:kitaev}) provide a straightforward
method to compute such examples. Here and in the following section 
we present compact formulas for the 8 examples obtained by
using Theorem \ref{lem:main}.

The same 8 examples (up to Okamoto and fractional-linear transformations)
are presented in \cite{Bo4}. Boalch uses transformation (\ref{qua:manin}),
so an Okamoto transformation has to be consequently applied to them in order
to get an algebraic function on the most motivating 
$\left(\frac12,\frac12,a,b\right)$ level. We strive to provide a detailed
and complementary account of application of quadratic transformations to
icosahedral Painlev\'e VI functions. In particular, we explain our ways of
representing algebraic Painlev\'e VI functions compactly.

Recall that a {\em Painlev\'e curve} is the normalization of an algebraic
curve defined by the minimal equation for an algebraic Painlev\'e VI solution
$y(t)$. The minimal equation is a polynomial in $y$ and $t$. The
indeterminant $t$ defines an algebraic map from the Painlev\'e curve to
$\PP^1$. As mentioned in \cite{Hi}, this map is a Belyi map. The reason is
that the corresponding field extension $\CC(t,y)\supset\CC(t)$ ramifies only
above $t=0,1,\infty$ due to the Painlev\'e property\footnote{Accordingly, in
\cite{Bo4} an algebraic Painlev\'e VI solution is defined as a triple
$(\Pi,y,t)$, where $\Pi$ is a compact curve, $y$ and $t$ are rational
functions on $\Pi$ such that $t$ is a Belyi function and $y(t)$ is a
Painlev\'e VI solution. However, this definition allows non-minimal Painlev\'e
curves, because Belyi maps can be appropriately composed to Belyi maps
again.}.

In this section we consider the transformations $39,40\Rightarrow 47,48$.
More transformations are presented in Appendix Section \ref{sec:other}.
Recall that Boalch classes which differ by replacing the local monodromy
differences $1/5$ with $2/5$ and vice versa are {\em siblings}. Such classes
are very similar. In particular, they have isomorphic Painlev\'e curves and
$t$-Belyi maps.

We start with the type-39 example in \cite[Section 7]{Bo2}, reparametrized
with $s\mapsto s-2$: 
\begin{eqnarray}
t_{39}\equal\frac12-\frac{(2s^7-18s^6+48s^5-50s^4+105s^3+3s^2-7s-3)u}{18(s^2-4s-1)(4s^2-s+1)^2},\\
y_{39}\equal
\frac12+\frac{14s^5-79s^4+6s^3+80s^2+116s-9}{6(s-1)(s^2-4s-1)\,u}.
\end{eqnarray}
Here $u=\sqrt{3(s+3)(4s^2-s+1)}$, so the function $y_{39}$ is defined on a
genus 1 curve. It a solution for
$P_{VI}\left(\frac13,\frac13,\frac45,\frac45; t_{39}\right)$.
Following Theorem \ref{lem:main}, we compute
\begin{equation}
\widetilde{y}_{39}=K_{\left[-\frac13,-\frac13,-\frac45,\,\frac45\,;\,t_{39}\right]}\,y_{39}=
\frac12-\frac{(2s^2-5s-1)\,u}{6\,(s-1)\,(4s^2-s+1)}.
\end{equation}
This is a solution of $P_{VI}\left(0,0,\frac7{15},\frac{13}{15};
t_{39}\right)$.  We compute:
\begin{eqnarray} %
\sqrt{t_{39}^{\,2}-t_{39}}\equal \frac{s\,(s+1)^2(s-2)^2(s-5)\sqrt{s\,(s+1)(s-2)(s+3)}}
{3\,(s^2-4s-1)(4s^2-s+1)\,u},\\
\sqrt{\widetilde{y}_{39}^{\;2}-\widetilde{y}_{39}}\equal
\frac{(s+1)\sqrt{s\,(s-2)(s+3)(s-5)}}{(s-1)\,u}.
\end{eqnarray}
We keep the factor $u$ in these expressions because we expect it will
disappear after simplifications. Note that the elliptic involution $u\mapsto
-u$ permutes the 2 singular points with zero local exponent differences.
This is the permutation which defines the subfield
$\CC(T_1,Y_1)\subset\CC(\tau,\eta)$ in terms of (\ref{eq:fields}).

The new square roots define the Painlev\'e curve for $y_{47}$. It appears to
be the fiber product of two elliptic curves:
\begin{equation}
v^2=s\,(s+1)(s-2)(s+3)\qquad \mbox{and}\qquad w^2=s\,(s-2)(s+3)(s-5).
\end{equation}
The fiber product is a hyperelliptic curve of genus 2. Its
Weierstrass form can be obtained by introducing the parameter
$q=\sqrt{(s-5)/(s+1)}$, so that $s=(q^2+5)/(1-q^2)$. Then the hyperelliptic
curve is represented by the equation
\begin{equation} \label{he47}
V^2=-(q^2+1)(q^2+5)(q^2-4).
\end{equation}
We can identify $v=6V/(q^2-1)^2$, $w=qv$. The connection with the
hyperelliptic form in \cite{Bo4} is $q=(j-3)/(j+3)$.

Formulas in (\ref{eq2:maint})--(\ref{eq2:mainsol}) gives us the following
solution of $P_{VI}\left(\frac12,\frac12,\frac13,\frac45;\,t_{47}\right)$:
\begin{eqnarray}
t_{47}\equal
\frac12-\frac{(2s^7-18s^6+48s^5-50s^4+105s^3+3s^2-7s-3)v}{4s^2(s+1)^3(s-2)^3(s-5)},\\
y_{47}\equal\nonumber
\frac12-\frac{(3s^7-27s^6+107s^5-205s^4+105s^3-37s^2-7s-3)v}
{2s(s+1)(s-2)^2(s-1)(3s^4-12s^3-14s^2-12s+3)}\\
&&-\frac{(4s^2-s+1)(s^2-4s-1)(7s^4-52s^3+34s^2-36s+15)}
{2s(s-2)^2(s-1)(3s^4-12s^3-14s^2-12s+3)\sqrt{(s+1)^3(s-5)}}.
\end{eqnarray}
Here we rewrote $y_{47}$ following the form (\ref{eq1:fmaint}).
This solution can be written in the variables $q$ and $V$ of the
hyperelliptic curve (\ref{he47}):
\begin{eqnarray*}
t_{47}\equal
\frac12-\frac{(q^{14}+17q^{12}+111q^{10}+370q^8+815q^6+1077q^4+169q^2+32)V}
{108q^2(q^2+1)^3(q^2+5)^2},\\
\label{sol:y47} y_{47}\equal
\frac12-\frac{(q^2-1)^2(q^4+7q^2+16)(q^4+7q^2+1)(q^8+8q^6+28q^4+36q^2-10)}
{6\,q(q^2+1)^2(q^2+2)(q^2+5)(q^8+8q^6+60q^4+176q^2-2)}\\
&&-\frac{(5q^{14}+85q^{12}+555q^{10}+1715q^8+2590q^6+1902q^4+1574q^2+322)V}
{6(q^2+1)^2(q^2+2)(q^2+5)(q^8+8q^6+60q^4+176q^2-2)}.
\end{eqnarray*}
Note that the fractional-linear transformation $\lambda\mapsto1-\lambda$ of
the isomonodromy problem (\ref{eq:JM}) acts on $y_{47}$ as the automorphism
$q\mapsto-q$, $V\mapsto-V$, which is not a hyperelliptic involution since it
has too few fixed points. Previously know examples of (hyper)elliptic
solutions of Painlev\'e VI have expressions on which the fractional-linear
transformation $\lambda\mapsto1-\lambda$ acts as a (hyper)elliptic
involution. However, there are type-47 solutions on which 
$\lambda\mapsto1-\lambda$ does act as a hyperelliptic involution.
For example, if we apply Theorem \ref{lem:main} to the fractional-linear transformation
$\pi_2\,y_{39}$ and use the same variables $q$, $V$,
then we get the following solution of
$P_{VI}\left(\frac12,\frac12,\frac15,\frac13;t_{47}\right)$:
\begin{gather*} \textstyle
\frac12+\frac{(3q^{14}+51q^{12}-6q^{11}+337q^{10}-84q^9+1168q^8-366q^7+2767q^6-534q^5+3821q^4-276q^3+413q^2-30q+80)V}{8q(q^2+1)(q^2+5)(q^{10}+10q^8-45q^7+25q^6-315q^5+35q^4-495q^3+80q^2-225q+11)}.
\end{gather*}

An example of type 48 can be computed similarly. With the same elliptic curve and
$t_{39}$ as just before, we have the following type-40 icosahedral function
(obtained after applying
$K_{[-\frac23,-\frac23,\frac25.\frac25;t_{39}]}\,\pi_2\,K_{[-\frac35,-\frac35,\frac23,\frac23;t_{39}]}$
and reparametrization $s\mapsto s-2$ to the corresponding example in
\cite{Bo2}):
\begin{equation}
y_{40}=\frac12+\frac{2s^6-14s^5+17s^4+16s^3-112s^2-2s-3}{2u\,(3s-1)(s^2-4s-1)}
\end{equation}
This is a solution of
$P_{VI}\left(\frac25,\frac25,\frac23,\frac23;t_{39}\right)$. Similar
application of Theorem \ref{lem:main} gives the following solution of
$P_{VI}\left(\frac12,\frac12,\frac25,\frac23;t_{47}\right)$.
\begin{eqnarray}  \label{sol:y48}
y_{48}\equal\nonumber
\frac12-\frac{(19s^6-138s^5+195s^4+380s^3+195s^2-138s-89)\,v}{2(s+1)^2(s-2)
(19s^5-155s^4+390s^3-590s^2-5s-3)}\\
&&\hspace{-8pt}+\frac{9(4s^2-s+1)(s^2-4s-1)(s^5-7s^4+13s^3-115s^2-2s-10)}
{2(s\!+\!1)^2(19s^5\!-\!155s^4\!+\!390s^3\!-\!590s^2\!-\!5s\!-\!3)
\sqrt{\!s(s\!-\!2)^3(s\!+\!3)(s\!-\!5)}}.
\end{eqnarray}
In terms of the variables $q$ and $V$, we have:
\begin{eqnarray*}
y_{48}\equal
\frac12+\frac{q^{15}+2q^{14}+17q^{13}+30q^{12}+103q^{11}+154q^{10}+242q^9+276q^8+K_0}
{12q(q^2+1)(4q^5+7q^4+28q^3+34q^2+34q+7)V},
\end{eqnarray*}
where $K_0=63q^7-806q^6-1327q^5-3322q^4-2295q^3-926q^2-260q-160$. On this
solution, the fractional-linear transformation $\lambda\mapsto 1-\lambda$ acts as a
hyperelliptic involution. But some other solutions of type 48 do not have this
``hyperelliptic" symmetry. As an example, one may the solution of
$P_{VI}\left(\frac12,\frac12,\frac13,\frac25\right)$ obtained after applying
Theorem \ref{lem:main} to $K_{[-\frac25-\frac25,\frac13,\frac13]}\pi_2\,y_{40}$.

\begin{figure}
\[ \begin{picture}(360,270)(0,16)
\put(24,20){$\PP^1_{39}$} \put(146,80){$\PP^1_\tau$} \put(23,220){$C_{39}$}
\put(265,20){$\PP^1_{47}$} \put(264,220){$C_{47}$} \put(146,280){$C_\tau$}
\put(30,81){\vector(0,-1){48}} \put(30,168){\line(0,1){46}}
\put(150,141){\vector(0,-1){48}} \put(150,229){\line(0,1){46}}
\put(270,81){\vector(0,-1){48}} \put(270,168){\line(0,1){46}}
\put(140,78){\vector(-2,-1){100}} \put(160,78){\vector(2,-1){100}}
\put(140,276){\vector(-2,-1){100}} \put(160,276){\vector(2,-1){100}}
\put(328,276){\vector(-1,-1){48}} \put(328,280){$C^*_{47}$}
\put(0,122){\fbox{ $\begin{array}{ccc} \!{\bf 5} & {\bf 5} & 5 \\
\!{\bf 5} & {\bf 5} & 3 \\ \!{\bf 3} & {\bf 3} & 3 \\
\!2 & 2 & 1,\!1\!\! \\ & & 1,\!1\!\! \end{array}$}}
\put(105,182){\fbox{ $\begin{array}{cccc} \!5 & 5 & 5 & 5 \\
\!5 & 5 & 3 & 3 \\ \!3 & 3 & 3 & 3 \\
\!1 & 1 & 1,\!1 & 1,\!1\!\! \\
\!1 & 1 & 1,\!1 & 1,\!1\!\! \end{array}$ }}
\put(235,122){\fbox{ $\begin{array}{ccc} 5 & {\bf 5} & {\bf 5}\! \\
5 & {\bf 3} & {\bf 3}\! \\
3 & {\bf 3} & {\bf 3}\! \\
\!1* & 2 & 2\! \\ \!1*& 2 & 2\! \end{array}$ }}
\end{picture} \]
\caption{Quadratic transformations for the 47th Boalch class}
\label{fig:sol47} \stepcounter{table}
\end{figure}
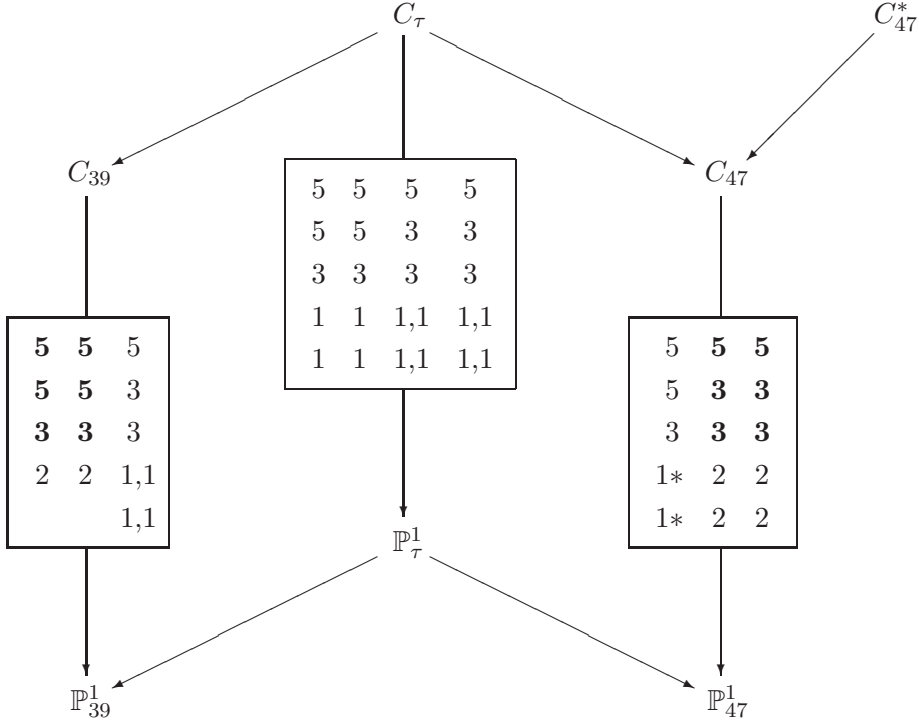
As mentioned above, the variables $t_{39}$ and $t_{47}$ define Belyi maps.
Figure \ref{fig:sol47} depicts change of branching of these Belyi maps. By $\PP^1_{39}$,
$\PP^1_{47}$ and $\PP^1_{\tau}$ we denote the projective lines with the rational parameters
$t_{39}$, $t_{47}$, $\sqrt{t_{39}}$, respectively. The curves $C_{39}$ and $C_{47}^*$ are
the Painlev\'e curves for, respectively, type 39 and 47 solutions. The function fields of
$C_{47}$ and $C_{\tau}$ are, respectively, $\CC(T,y_{39})=\CC(s,v)$ and
$\CC(t,T,y_{39})=\CC(s,u,v)$. The
map $t_{39}$ is from $C_{39}$ to $\PP^1_{39}$. The map $t_{47}$ is from
$C^*_{47}$ to $\PP^1_{47}$. Non-vertical arrows represent degree 2
coverings. In boxes we represent the branching patterns of the morphisms.
Each column gives branching orders of one fiber. In the middle box, the
first two columns represent points with $s\in\{-1,2,0,5\}$; the last two
columns represent the points with $s=\infty$, $4s^2=s+1$ or $s^2=4s+1$. The
bold numbers represent branching points of the upper degree 2 coverings
above the parallelograms. The stars mark the branching points of the
upper-right degree 2 covering. Each parallelogram is a commutative fiber
product diagram. The genus of $C_\tau$ is 4. The two composite coverings
have the following branching pattern, respectively:
\[
\begin{array}{ccc} 10 & 10 & 5,\!5 \\
10 & 10 & 3,\!3 \\ 6 & 6 & 3,\!3 \\
\!2 & 2 & 1,\!1,\!1,\!1 \\ 2 & 2 & 1,\!1,\!1,\!1
\end{array}\qquad\mbox{and}\qquad
\begin{array}{ccc} 5,\!5  & 10 & 10 \\
5,\!5  & 6 & 6 \\ 3,\!3  & 6 & 6 \\
1,\!1  & 2,\!2 & 2,\!2 \\ 1,\!1 & 2,2 & 2,2
\end{array}.
\]
The quadratic covering $\PP^1_{\tau}\to\PP^1_{39}$ branches above the two
points represented by the first two columns of the first box. The covering
$\PP^1_{\tau}\to\PP^1_{47}$ branches above the two points represented by the
last two columns of the third box.

\section{More algebraic Painlev\'e VI functions}
\label{sec:morealg}

Here we present compact expressions for the remaining 6 difficult
cases of icosahedral Painlev\'e VI functions. We have the following cascades of
quadratic transformations:
\begin{equation}
31,32\Rightarrow 44,45\Rightarrow 50,51;\qquad  41\Rightarrow 49\Rightarrow
52.
\end{equation}
They start with an icosahedral case Okamoto equivalent to a
Dubrovin-Mazzocco example in \cite{DM}, and end with an icosahedral case
$\left(\frac12,\frac12,\frac12,b\,\right)$. This is not surprising since
\cite{DM} classifies all algebraic solutions of the equations
$P_{VI}(0,0,0,B;t)$, to which quadratic transformation (\ref{qua:manin}) can
be applied consequently twice.

First we recall expressions for the Dubrovin-Mazzocco examples transformed to the
$(a,a,b,b)$ level. A solution of type-31 equation
$P_{VI}\left(\frac15,\frac15,\frac15,\frac15;t_{31}\right)$ is
\begin{equation}
t_{31}=\frac{(s+1)^5(s-3)^3(s^2+4s-1)}{256s^3(s^2-5)},\qquad
y_{31}=-\frac{(s+1)^4(s-3)^2}{4s(s^2+3)(s^2-5)}.
\end{equation}
This is reparametrization $s\mapsto (s-3)/(s+1)$ of an expression in
\cite[Section 3.4]{K3}. A solution of type-32 equation
$P_{VI}\left(\frac25,\frac25,\frac25,\frac25;t_{32}\right)$ is given by
\begin{equation}
t_{32}=t_{31},\qquad\qquad
y_{32}=-\frac{(s+1)^2(s-3)^2(s^2+4s+7)}{48s(s^2-5)}.
\end{equation}
This is reparametrization $s\mapsto (1-s)/(s+3)$ of
$K_{[-\frac25,-\frac25,-\frac25,-\frac25]}K_{[0,0,0,-\frac45]}H_3$, where $H_3$
is the Icosahedron solution presented in \cite{DM}. Finally, a solution of type-41 equation
$P_{VI}\left(\frac13,\frac13,\frac13,\frac13;t_{41}\right)$ is parametrized in
\cite[Theorem C]{Bo2}:
\begin{eqnarray}
t_{41}\equal\textstyle
\frac12+\frac{(s+1)(32s^8-320s^7+1112s^6-2420s^5+3167s^4-2420s^3+1112s^2-320s+32)}
{54u^3s(s-1)},\\
y_{41}\equal\textstyle
\frac12-\frac{8s^7-28s^6+75s^5+31s^4-269s^3+318s^2-166s+56}{18u(s-1)(3s^3-4s^2+4s+2)},
\end{eqnarray} where $u=\sqrt{s(8s^2-11s+8)}$.

Application of Theorem \ref{lem:main} to the functions $y_{31}$ and $y_{32}$
gives algebraic expressions for icosahedral functions
of types 44, 45  defined over the field
\begin{equation}
\CC\left(s^2,\sqrt{(s^2-1)(s^2-9)},\sqrt{s^4-18s^2+1}\right).
\end{equation}
The invariance under $s\mapsto -s$ is explained by the fact that this is the involution
$\tau\to1/\tau$, $\eta\to1/\eta$ in terms of the diagrams as (\ref{eq:fields});
in particular, it realizes the fractional-linear transformation $\lambda\to1-\lambda$.
Hence $s\mapsto -s$ must indeed fix the function fields $\CC(t_{44},y_{44})$, $\CC(t_{45},y_{45})$.
With $s^2$ effectively being a rational parameter for the transformed functions,
the square root $\sqrt{s^4-18s^2+1}$ defines a genus 0 curve.
We can parametrize it by setting \mbox{$s^2\mapsto q(2q+1)/(q-2)$}. Then we can identify:
\begin{eqnarray*}
\sqrt{s^4-18s^2+1}=\frac{2(q^2-4q-1)}{q-2},\quad
\sqrt{(s^2-1)(s^2-9)}=\frac{2\sqrt{(q^2+1)(q^2-4q+9)}}{q-2}.
\end{eqnarray*}
Evidently, the transformed solutions of the types 44 and 45 are defined on the elliptic curve
$u^2=(q^2+1)(q^2-4q+9)$. Compared with the respective parametrization in \cite{Bo4}, we
have $s\mapsto -1/s$ and $j\mapsto 1-4/q$. Eventually, the functions
$y_{31}$ or $y_{32}$ are transformed to the following solutions of
$P_{VI}\left(\frac12,\frac12,\frac15,\frac15;t_{44}\right)$ and
$P_{VI}\left(\frac12,\frac12,\frac25,\frac25;t_{44}\right)$:
\begin{eqnarray} \label{sol:t44}
t_{44}\equal\textstyle\frac12+\frac
{q^{10}-10q^9+45q^8-120q^7+190q^6+4q^5-410q^4+680q^3+25q^2-90q-27}
{2(q^2-4q-1)(q^2+1)u^3},\\ 
y_{44}\equal\textstyle\frac12+\frac
{7q^{10}-70q^9+329q^8-908q^7+1494q^6+24q^5-3310q^4+5692q^3+211q^2-418q-171}
{2(q^2-4q-1)(7q^6-28q^5+91q^4-88q^3+229q^2-76q+57)u},\\
\label{sol:y45} y_{45}\equal\textstyle\frac12
+\frac{3q^{10}-24q^9+93q^8-216q^7+446q^6+320q^5-1342q^4+2824q^3+527q^2-312q-207}
{6(q^2-4q-1)(q^4-2q^3+8q^2+10q+23)(q^2+1)u}.
\end{eqnarray}
Surely, we may choose to parametrize the square root $\sqrt{(s^2-1)(s^2-9)}$ instead,
say $s^2\mapsto (q^2-9)/(q^2-1)$.  But then the transformation $\lambda\to1-\lambda$
would not be acting on the obtained expressions for $y_{44}$ and $y_{45}$ as
flipping the sign of the other root $\sqrt{s^4-18s^2+1}$, so the ``elliptic" symmetry
would be lost. However, if we would consider the
quadratic transformation of $\pi_2y_{31}$ or $\pi_2y_{32}$, we would obtain
a solution of $P_{VI}\left(\frac12,\frac12,\frac15,\frac45;t_{44}\right)$ or
$P_{VI}\left(\frac12,\frac12,\frac25,\frac35;t_{44}\right)$, on which $\lambda\to1-\lambda$
acts by flipping the sign of $\sqrt{s^4-18s^2+1}$ but not of $\sqrt{(s^2-1)(s^2-9)}$;
then we should use the latter parametrization of $s^2$ in order to obtain expressions with
the ``elliptic" symmetry.

Application of Theorem \ref{lem:main} to the functions $y_{44}$, $y_{45}$
gives us expressions for icosahedral functions of types 50, 51 defined over the field
\begin{equation}
\CC\left(q,\sqrt{q(q-2)(2q+1)},\sqrt{-(2q+1)(q^2-2q+5)}\right).
\end{equation}
For comparison with the previous case,
notice that the square root $u$ present in (\ref{sol:t44})--(\ref{sol:y45}) is ``killed" by the
the involution $\tau\mapsto1/\tau$, $\eta\mapsto1/\eta$ of (\ref{eq:fields}). The function
field defines an algebraic curve of genus 3; the curve is not hyperelliptic.
As found in \cite{Bo4}, the Painlev\'e curve can be defined by the equation
\begin{equation} \label{boalc50}
5p^4+6p^2r^2+5r^4+6p^2+6r^2+1=0.
\end{equation}
With reference to (\ref{sol:t44})--(\ref{sol:y45}), we offer the following ``parametrization":
\begin{equation} \label{subs50}
p=\sqrt{-\frac{2q+1}{q^2-2q+5}},\qquad r=\sqrt{-\frac{q(q-2)}{q^2-2q+5}}.
\end{equation}
The variable $q$ can be expressed as
\begin{equation}
q=\frac{5p^2+5r^2+1}{2p^2}\qquad\mbox{or}\qquad
q=-\frac{5p^2+r^2+1}{2(r^2+1)}.
\end{equation}
It appears that application of Theorem \ref{lem:main}
to $\pi_2y_{44}$ and $\pi_2y_{45}$ 
gives simpler icosahedral functions of type 50, 51.
Here are expressions for those solutions, of
$P_{VI}\left(\frac12,\frac12,\frac12,\frac45;t_{50}\right)$ and
$P_{VI}\left(\frac12,\frac12,\frac12,\frac35;t_{50}\right)$:
\begin{eqnarray} \label{t50pqr}
t_{50}\equal\textstyle\frac12+
\frac{32p^6+32p^4r^2+56p^4+64p^2r^2-4p^2-68r^2+q^4-4q^3-6q^2+20q-3}
{16p^3rq(q-2)^2},\quad\\
y_{50}\equal\textstyle\frac12+\frac{p\left(rK_3-40p^4-200p^2r^2+355r^2+261p^2
+2q^5-6q^4-15q^3+23q^2-33q+52\right)
}{q(q-2)(2q+1)(3p^2+r^2+2r+1)(5r^3+7p^2r-3p^2q+11p^2+2r)},\\
\label{eq:y51} y_{51}\equal\textstyle\frac12-
\frac{120p^5+200p^3r^2-120p^4-200p^2r^2+132p^3+460pr^2-180p^2-100r^2
+K_4}{2r(2q+1)(24p^3+40pr^2-48p^2-80r^2+2pq^2+12pq+2p-7q^2-2q-11)},
\end{eqnarray}
where $K_3=440p^4+200p^2r^2+952p^2+760r^2+4q^5-22q^4+15q^3+66q^2+89q+174$,
and $K_4=p(q^5-6q^4+7q^3+12q^2-126q+24)-q^5+6q^4-4q^3-24q^2-33q-30$.
To have most compact expressions for these functions, we use in (\ref{t50pqr})--(\ref{eq:y51})
all three variables $q,p,r$.
To find these expressions, we followed this strategy: first we factored over $\QQ$
the divisors on the Painlev\'e curve defined by the numerators and denominators
of a rational expression under consideration;
for each divisor irreducible over $\QQ$ we found a Gr\"{o}bner
basis of polynomial functions vanishing on it; then we tried to combine
lowest degree polynomials from the Gr\"{o}bner bases of divisor factors so to
built a compact expression for the same rational function; and finally we tried to
change a monomial basis for each chosen polynomial so to reduce the size of coefficients.
To illustrate our approach, we present another expression for $t_{50}$, in $p$ and $r$ only:
\begin{equation}\label{t50pr} \textstyle
t_{50}=\frac12+\frac{32r^6p^4-36r^6p^2-36r^2p^6+32r^4p^6-38r^6-30r^4p^2-30r^2p^4
-38p^6-13r^4+6p^2r^2-13p^4}{40pr(p^2+1)(r^2+1)(5p^2+5r^2+1)}.
\end{equation}
Notice that the large numerators in (\ref{t50pqr}) and (\ref{t50pr}) both contain 11 terms:
it is typical that changing monomial basis does not change the number of terms in
a polynomial expression. Even in the expression of $t_{52}$ in $q$ only, easily
obtainable by substituting (\ref{subs50}), the numerator has degree 10 in $q$, hence 11 terms.
The numerator of (\ref{t50pqr})
was obtained after elimination of the monomial $q^5$
from two polynomials in a Gr\"obner basis with respect to a total degree ordering.
The monomial basis in the numerator of (\ref{t50pr}) was intentionally chosen symmetric.
Incidentally, the whole expression in (\ref{t50pr}) is symmetric in $p,r$.
For better illustration of intermediate steps of our strategy, recall  that $t_{50}$ must be a Belyi map,
so it must be highly factorizable. 
The divisor of $t_{50}$ as a rational function on the curve (\ref{boalc50}) can be found to be
\begin{equation}   \label{divexpr}
\begin{smallmatrix}
5\,(p+r-1,\,2r^2-2r+1)+5\,(p+r+1,\,2r^2+2r+1)+3\,(p+r-1,\,2r^2-2r+3)+3\,(p+r+1,\,2r^2+2r+3)\\
+2\,(p-r,\,16r^4+12r^2+1)-5\,(r,\,p^2+1)-5\,(p,\,r^2+1)-3\,(r,\,5p^2+1)-3\,(p,\,5r^2+1)
-2\,(5p^4+6p^2r^2+5r^4,\,\infty).
\end{smallmatrix}
\end{equation}
Here we represented each irreducible divisor by a Gr\"obner basis of it, and we used the symbol
$\infty$ to indicate the points at infinity. The multiplicity pattern is consistent with the respective
entry in the last column of Table 1 in \cite{Bo2}. Eventually, by combining appropriate
polynomial factors (and adjusting a scalar multiple) in the same manner as in
\cite[Section 4]{Vid}, we can arrive at the following compact expression:
\begin{equation}
t_{50}=\frac{5(p-r)^2(p+r-1)^3(p+r+1)^3(3p^2+3r^2-2pr+1)}
{128pr(p^2+1)(r^2+1)(5p^2+5r^2+1)}.
\end{equation}
In general, this procedure is tricky and hard to automatize in an effective way.

Application of Theorem \ref{lem:main} to the function $y_{41}$
gives us an icosahedral function of type 49 defined over the field
\begin{equation}
\CC\left(s,\sqrt{(s-2)(2s-1)(2s^2+s+2)},\sqrt{s^2-7s+1}\right).
\end{equation}
Since we can parametrize the latter square root, the Painlev\'e curve is hyperelliptic.
We choose the following parametrization by $q$:
\begin{equation}
s\mapsto -\frac{(\widetilde{q}-1)(\widetilde{q}-5)}{(\widetilde{q}+1)(\widetilde{q}+5)}, \qquad
\widetilde{q}\mapsto \sqrt{5}\,q.
\end{equation}
The hyperelliptic curve (of genus 3) is $v^2=3(q^4+14q^2+1)(5q^4+6q^2+5)$.
The relation with the parameter in \cite{Bo4} is
$j=2(\widetilde{q}-2)/(\widetilde{q}+1)$. We get the following solution of
$P_{VI}\left(\frac12,\frac12,\frac13,\frac13;t_{49}\right)$:
\begin{eqnarray} \label{sol:t49}
t_{49}\equal\textstyle\frac12+
\frac{q\,(1125q^{16}+7000q^{14}+26124q^{12}+11112q^{10}-25186q^8+11112q^6+26124q^4+7000q^2+1125)}
{(q^2-1)(5q^4+6q^2+5)\,v^3},\quad\\
\label{sol:y49} y_{49}\equal\textstyle\frac12+
\frac{8qK_1+\sqrt5\,(q^2-1)(q^2+1)(3q^2+1)(q^2+3)(q^4+14q^2+1)(5q^6+53q^4-9q^2+15)}
{4v(q^2-1)\,\left(K_2\,+\,8\sqrt5\,q\,(q^2-1)(q^2+1)^2(3q^2+1)(q^2+3)\right)},
\end{eqnarray}
where $K_1=
200q^{16}\!+\!1367q^{14}\!+\!4835q^{12}\!+\!3643q^{10}\!-\!1609q^8\!+
\!2933q^6\!+\!4025q^4\!+\!825q^2\!+\!165$,
and $K_2=35q^{12}+450q^{10}+1097q^8+1324q^6+805q^4+370q^2+15$.

Application of Theorem \ref{lem:main} to $y_{49}$ or $\pi_2y_{49}$ gives us
icosahedral functions of type 52 defined over the field
\begin{equation}
\CC\left(q,\sqrt{-2(q^2+1)(5q^2-1)(q^2+3)},\sqrt{2(q^2+1)(q^2-5)(3q^2+1)}\right).
\end{equation}
The function
field defines an algebraic curve of genus 7; the curve is not hyperelliptic.
As found in \cite{Bo4}, the Painlev\'e curve can be defined by the equation
\begin{eqnarray}
9p^6r^2+18p^4r^4+9p^2r^6+4p^6+26p^4r^2+26p^2r^4+4r^6+\nonumber\\
8p^4+57p^2r^2+8q^4+20p^2+20r^2+16 \equal 0.
\end{eqnarray}
With reference to (\ref{sol:t49})--(\ref{sol:y49}), we offer the following ``parametrization"
\begin{equation}
p=\sqrt\frac{(q+\sqrt5)(3q^2+1)}{2(q-\sqrt5)(q^2+1)},\qquad
r=\sqrt\frac{-(\sqrt5q-1)(q^2+3)}{2(\sqrt5q+1)(q^2+1)}.
\end{equation}
In the other direction, we have
\begin{equation}
q=\frac{5p^4+6p^2r^2+r^4+11p^2-r^2+8}{\sqrt5\,(p^2-r^2)(p^2+r^2-1)}.
\end{equation}
We choose to present the quadratic transformation of $\pi_2y_{49}$, since it looks simpler
than quadratic transformation of $y_{49}$. Here is the solution of
$P_{VI}\left(\frac12,\frac12,\frac12,\frac23;t_{52}\right)$:
\begin{align*}
t_{52}\;
\equal&\textstyle\!\frac12-\frac{(p^2+r^2-1)
(p^8-7p^6r^2-24p^4r^4-7p^2r^6+r^8+6p^4r^2+6p^2r^4+14p^4+p^2r^2+14r^4+16p^2+16r^2+17)}
{96p^3r^3\left(p^2+r^2+1\right)},\\
y_{52}\;\equal& \textstyle \frac12+
\frac{K_5-4800p^6-28160p^4+8704q^2r^2-6240q^4-60352p^2-9736q^2-7288
-5r(966q^4-512r^2+608q^2-982)}
{4\sqrt5p\,r^2(q+\sqrt5)(\sqrt5q+1)^2\left(K_6-60p^2-48q^2+22-20r(p^2q^2+3p^2+4)\right)},
\end{align*}
where
{\footnotesize
\begin{eqnarray*}
\sqrt{5}qK_5\!\!&\!\!=\!\!&\!\!
96q^2r^6\!-17256p^4q^4\!+288r^6\!-4992p^4r^2\!-14400p^6\!-2175q^6\!-5504p^2r^2\!-1752p^4\\
&&\!\!-39859q^4\!+51136p^2\!-6805q^2\!+25703-5r(3360p^6q^2\!+5312p^4q^4\!-1024r^4q^4\!+6240p^6\\
&&\!\!+12768p^4q^2\!-285q^6\!+5600p^4\!-2088r^2q^2\!+8319q^4\!-504r^2\!-4819q^2\!+1),\\
\sqrt{5}qK_6\!\!&\!\!=\!\!&\!\!
228q^2r^4\!-16r^2q^2\!-228p^2r^2\!-148p^2\!-56r^2\!-369q^2\!-281
-20r(3r^2q^2\!+12p^2\!+5r^2\!+11q^2\!+9).
\end{eqnarray*}}
\hspace{-9pt}
There might be shorter expressions for $y_{52}$  with other monomial bases for the
large polynomials in the numerator and denominator, but it combinatorially hard
to find such monomial bases.

\section{Appendix}
\label{sec:other}

First we present Table \ref{fig:frli} of the fractional-linear
transformations for Painlev\'e VI functions. If one starts with a solution
$y(t)$ of $P_{VI}(a,b,c,d+1)$, in each row we give a solution of a
Painlev\'e VI equation with permuted singular points
(of the corresponding Fuchsian system) in terms of $y(t)$ and $t$.
We also give expressions of the transformations in \cite{TOS} notation,
and fractional-linear transformations of $\lambda$ for 
Fuchsian system (\ref{eq:JM}). 
There might be non-trivial Schlessinger transformations of $\Psi$ as well
(if $\lambda=\infty$ is moved). The substitutions
of $\lambda$ compose as a direct group action --- from right to left.
The $\theta$'s, $y$ and $t$ are functions on the singular points,
so the substitutions for them compose from left to right.
As mentioned, we use the notation $\pi_2$ in Sections \ref{sec:algebraic} and \ref{sec:morealg}.
Note that non-trivial transformations which fix the argument $t$ correspond
to permutations of the conjugacy class $2+2$. We actually use only these
transformations, or change the argument to $t/(t-1)$.
\begin{table}\begin{center}
\begin{tabular}{|c|c|c|c|c|c|c|} \hline
\cite{TOS} notation & $\lambda$ &
$(\theta_0,\theta_1,\theta_t,1\!-\!\theta_{\infty})$
 & $y$ & $t$ \\ \hline
id & $\lambda$ & $(a,b,c,d)$ & $y$ & $t$ \\
$\sigma_1\pi_2$ & $t\lambda/(\lambda+t-1)$ & $(a,b,d,c)$ & $(1-t)y/(y-t)$ & $1-t$ \\
$\sigma_2\pi_1$ & $t\lambda$ & $(a,c,b,d)$ & $y/t$ & $1/t$ \\
$\sigma_1\sigma_2\pi_2$ & $t\lambda/(t\lambda+1-t)$ & $(a,c,d,b)$ & $(t-1)y/t(y-1)$ & $(t-1)/t$ \\
$\sigma_2\sigma_1\pi_1$ & $t\lambda/(\lambda-1)$ & $(a,d,b,c)$ & $y/(y-t)$ & $1/(1-t)$ \\
$\sigma_1\sigma_2\sigma_1$ & $\lambda/(\lambda-1)$ & $(a,d,c,b)$ & $y/(y-1)$ & $t/(t-1)$ \\
$\sigma_1$ & $1-\lambda$ & $(b,a,c,d)$ & $1-y$ & $1-t$ \\
$\pi_2$ & $t(\lambda-1)/(\lambda-t)$ & $(b,a,d,c)$ & $t(y-1)/(y-t)$ & $t$ \\
$\sigma_1\sigma_3$ & $t\lambda-\lambda+1$ & $(b,c,a,d)$ & $(y-1)/(t-1)$ & $1/(1-t)$ \\
$\sigma_3\pi_1$ & $t/(\lambda-t\lambda+t)$ & $(b,c,d,a)$ & $t(y-1)/(t-1)y$ & $t/(t-1)$ \\
$\pi_2\sigma_2$ & $(t\lambda-1)/(\lambda-1)$ & $(b,d,a,c)$ & $(y-1)/(y-t)$ & $1/t$ \\
$\sigma_1\sigma_2$ & $1/(1-\lambda)$ & $(b,d,c,a)$ & $(y-1)/y$ & $(t-1)/t$ \\
$\sigma_3\sigma_1$ & $t(1-\lambda)$ & $(c,a,b,d)$ & $(t-y)/t$ & $(t-1)/t$ \\
$\sigma_2\pi_2$ & $t(\lambda-1)/(t\lambda-1)$ & $(c,a,d,b)$ & $(y-t)/t(y-1)$ & $1/t$ \\
$\sigma_3$ & $\lambda-t\lambda+t$ & $(c,b,a,d)$ & $(y-t)/(1-t)$ & $t/(t-1)$ \\
$\sigma_3\sigma_2$ & $t/(t\lambda-\lambda+1)$ & $(c,b,d,a)$ & $(y-t)/(1-t)y$ & $1/(1-t)$ \\
$\pi_1\pi_2$ & $(\lambda-t)/(\lambda-1)$ & $(c,d,a,b)$ & $(y-t)/(y-1)$ & $t$ \\
$\sigma_1\pi_1$ & $t/(1-\lambda)$ & $(c,d,b,a)$ & $(y-t)/y$ & $1-t$ \\
$\sigma_3\pi_2$ & $t(\lambda-1)/\lambda$ & $(d,a,b,c)$ & $t/(t-y)$ & $t/(t-1)$ \\
$\sigma_2\sigma_1$ & $(\lambda-1)/\lambda$ & $(d,a,c,b)$ & $1/(1-y)$ & $1/(1-t)$ \\
$\sigma_2\sigma_3$ & $(t\lambda-t+1)/\lambda$ & $(d,b,a,c)$ & $(1-t)/(y-t)$ & $(t-1)/t$ \\
$\sigma_2$ & $1/\lambda$ & $(d,b,c,a)$ & $1/y$ & $1/t$ \\
$\pi_1\sigma_1$ & $(\lambda+t-1)/\lambda$ & $(d,c,a,b)$ & $(t-1)/(y-1)$ & $1-t$ \\
$\pi_1$ & $t/\lambda$ & $(d,c,b,a)$ & $t/y$ & $t$ \\
\hline
\end{tabular}\end{center}
\vspace{-8pt}\caption{Fractional-linear transformations}\label{fig:frli}
\end{table}

Further in this Appendix, we present:
\begin{itemize}
\item The relation of our notation to the notation in \cite[Section 3]{TOS}. {(Compared with
\cite{Bo2} and \cite{Bo4}, the parameters $\theta_1,\theta_2,\theta_3,\theta_4$ there
should be identified with our $\theta_0,\theta_t,\theta_1,\theta_\infty$, respectively.)}
\item Lemma \ref{lem:schlesinger} on the relation between Okamoto
transformations and integer shifts in the local monodromy differences.
 \item Alternative notation to work with the quadratic transformations of
Painlev\'e VI functions.
\end{itemize}

The notation in \cite[Section 3]{TOS} can be related to our
notation as follows:
\begin{equation}
(\theta_0,\,\theta_1,\,\theta_t,\,\theta_\infty)\longleftrightarrow(-\alpha_4,-\alpha_3,-\alpha_0,1-\alpha_1).
\end{equation}
Then Okamoto transformation (\ref{eq:oka}) corresponds to the transformation
$s_2$ of \cite{TOS}. The transformations $s_0,s_1,s_3,s_4$ change the sign
of, respectively, $\theta_t,1-\theta_\infty,\theta_1,\theta_0$. In our 
paper, the action of these four transformations on the Painlev\'e VI
equation is considered trivial. 
The quadratic transformation $\psi_{VI}^{[2]}$ of \cite{TOS} converts
$\left(\frac{C-1}{2},\frac{B}2,\frac{B}2,\frac{C+1}2\right)$ to $(0,0,B,C)$,
which is the opposite direction 
from the accustomed here.
The quadratic transformation $(\frac12,a,b,\frac12)\mapsto(a,a,b,b)$ can be
realized by one of the following compositions:
\begin{equation} \label{eq:tosnot}
s_0\,s_2\,\psi_{VI}^{[2]}\,s_3\,s_2\,s_0\,s_4,\qquad \mbox{or}\qquad
s_0\,s_3\,s_4\,s_2\,\psi_{VI}^{[2]}\,s_3\,s_2\,s_0\,s_3\,s_4.
\end{equation}
Transformations (\ref{qua:manin0}) and (\ref{qua:manin}) can be realized
as $(\psi_{VI}^{[2]}\,\sigma_1)^{-1}$ and $(\sigma_1\,\sigma_2\,\psi_{VI}^{[2]})^{-1}$,
respectively.
In the notation of \cite{TOS}, our main formulas can be viewed (up to
fractional-linear transformations) as the result of eliminating
differentiation from (\ref{eq:tosnot})
and $\psi_{VI}^{[2]}s_3s_2s_0s_4$. As mentioned in Introduction,
we try to avoid composing differentiation and quadratic transformations
in proofs of our results, which leads us to composing (\ref{qua:manin}) with
contiguous relations within two Okamoto orbits.

The following lemma is basically noticed in \cite{O}. We present the observation in
the most direct notation. Recall that pairs of Fuchsian systems (\ref{eq:JM})
whose local monodromy differences differ by integers shifts as in this lemma,
with $k_0+k_1+k_t+k_{\infty}$ even, are related by
Schlesinger transformations \cite{JM}.
\begin{lemma} \label{lem:schlesinger}
Suppose that $k_0,k_1,k_t,k_{\infty}$ are integers. Then two Painlev\'e VI
equations $P_{VI}\left(\nu_0,\nu_1,\nu_t,\nu_{\infty};t\right)$ and
$P_{VI}\left(\nu_0\!+\!k_0,\nu_1\!+\!k_1,\nu_t\!+\!k_t,\nu_{\infty}\!+\!k_\infty;t\right)$
are related by a chain of Okamoto transformations if and only if either all
four integers are even, or they are all odd. The two Painlev\'e VI equations
are related by a combination of Okamoto and fractional-linear
transformations if and only if the sum $k_0+k_1+k_t+k_{\infty}$ is even.
\end{lemma}
\begin{proof} Let us denote $\Theta=(\nu_0+\nu_1+\nu_t+\nu_\infty)/2$. The composition
$$K_{\left[\Theta-\nu_0,\,\Theta-\nu_1,\,\Theta-\nu_t,\,2+\Theta-\nu_\infty;\,t\right]}\,
K_{\left[\nu_0,\;\nu_1,\;\nu_t,\;\nu_\infty;\;t\right]}$$ increases all four
local monodromy differences by 1. If we change the sign of $\nu_0$, $\nu_1$
or $\nu_t$ (in the definition of $\Theta$ as well), we change the shift sign
in the corresponding local monodromy difference. To have a negative shift at
$\infty$, one may consider the inverse transformations. The ``if" part of
the first statement follows. Suppose that $k'_0,k'_1,k'_t,k'_{\infty}$ are
integers, either all even or all odd. Okamoto transformations of
$P_{VI}\left(\nu_0\!+\!k'_0,\nu_1\!+\!k'_1,\nu_t\!+\!k'_t,\nu_{\infty}\!+\!k'_\infty;t\right)$
have the following forms, up to parameter permutations
(but not fractional-linear permutations!)
of the conjugacy class $2+2$:
\begin{eqnarray*}
&&\textstyle P_{VI}\!\left( \Theta-\nu_0+k''_0,\, \Theta-\nu_1+k''_1,\,
\Theta-\nu_t+k''_t,\, \Theta-\nu_\infty+k''_{\infty};\, t \right) \quad\mbox{or}\\
&&\textstyle P_{VI}\!\left( \Theta+k''_0, \Theta\!-\!\nu_0\!-\!\nu_1+k''_1,
\Theta\!-\!\nu_0\!-\!\nu_t+k''_t, \Theta\!-\!\nu_1\!-\!\nu_t+k''_{\infty};\,
t \right),
\end{eqnarray*}
with the same restrictions on $k''_0,k''_1,k''_t,k''_{\infty}$. Okamoto
transformations of these two equations keep the form of all three equations,
up to parameter permutations of the conjugacy class $2+2$. The ``only if"
part follows as well.

The second statement follows from the first one, since a shift by two odd
and two even integers can be obtained by composing a shift by all odd or all
even integers with a fractional-linear transformation which moves $\infty$.
For example \cite[Theorem 6]{GF}, a simplest Schlesinger transformation
$\left(\nu_0,\nu_1,\nu_t,\nu_{\infty}\right)\mapsto
\left(\nu_0+1,\nu_1,\nu_t,\nu_{\infty}+1\right)$ can be realized as
$$L\,K_{\left[\nu_0-\Theta,\,\Theta-\nu_1,\,\Theta-\nu_t,\,\nu_\infty-\Theta;\,t\right]}\,
K_{\left[\nu_0,\;\nu_1,\;\nu_t,\;\nu_\infty;\;t\right]},$$
 where $L$ is the fractional-linear transformation $\lambda\mapsto
t/\lambda$.
\end{proof}

Now we consider the quadratic transformations in an alternative notation. First we concentrate
on quadratic transformations within Table \ref{fig:quadr}. 
Consider $t_2$,  $g_1$, $T_2$ as in the context of Theorem \ref{lem:main};
we relate to $\eta,\tau$ and other ``global" notation via (\ref{eq:frlis})--(\ref{eq:tauet}).
Let $G_3$ denote the solution $(\tau+1)(\eta-1)/2(\eta-\tau)$
of $P_{VI}(\frac{A}2,\frac{A}2,\frac{B}2,\frac{B}2+1;T_2)$. We express the functions and
variables as follows:
\begin{equation}
g_1=\frac{1+\psi}2, \qquad t_2=\frac{1+\theta}2, \qquad G_3=\frac{1+\xi}2,
\qquad T_2=\frac{1+\sigma}2.
\end{equation}
Like we commented before Theorems \ref{the:main} and \ref{lem:main},
the points $\lambda=0$ and $\lambda=1$ are symmetric in this form, since the fractional-linear
permutation of them merely changes the sign of $\psi,\theta$ or $\xi,\sigma$. Hence this form
is widely used in \cite{K3}, \cite{Bo2} and Sections \ref{sec:algebraic} and \ref{sec:morealg}
here. (The unindexed $\theta$ does not denote any local monodromy differences from here on.)
In terms of $\tau$ and $\eta$, we have:
\begin{equation} \label{eq:newvars}
\theta=\frac{\tau^2+1}{\tau^2-1},\qquad
\psi=\frac{\eta^2+1}{\eta^2-1},\qquad
\sigma=\frac12\!\left(\tau+\frac1{\tau}\right),\qquad
\xi=\frac{\tau\eta-1}{\eta-\tau}.
\end{equation}
We can identify
\begin{eqnarray} \label{eq:sqid1}
\sqrt{\theta^2-1}=\frac{2\tau}{\tau^2-1},\qquad
\sqrt{\psi^2-1}=\frac{2\eta}{1-\eta^2},\qquad
\sqrt{\sigma^2-1}=\frac12\!\left(\tau-\frac1{\tau}\right),
\end{eqnarray}
so that $\sqrt{t_2^2-t_2}=\frac12\sqrt{\theta^2-1}$ and
$\sqrt{g_1^2-g_1}=\frac12\sqrt{\psi^2-1}$ in the setting of (\ref{eq:tauet}).
In particular, reminiscent to Corollary 3 in \cite{Bo4}, we have
\begin{equation}
\sigma=\frac{\theta}{\sqrt{\theta^2-1}},\qquad\qquad
\xi=\frac{\sqrt{\psi^2-1}-\sqrt{\theta^2-1}}{\psi-\theta}.
\end{equation}
To transform from $G_3$ to $g_1$, the formulas are:
\begin{equation}
\theta=\frac{\sigma}{\sqrt{\sigma^2-1}},\qquad\qquad
\psi=\frac{\sigma\xi^2-2\xi+\sigma}{(\xi^2-1)\,\sqrt{\sigma^2-1}}.
\end{equation}
We also have the following relations:
\begin{eqnarray}
\sqrt{\theta^2-1}=\frac{1}{\sqrt{\sigma^2-1}},\hspace{15pt}&&
\sqrt{\psi^2-1}=\frac{\xi^2-2\sigma\xi+1}{(1-\xi^2)\,\sqrt{\sigma^2-1}},\\
\label{eq:sqid2}
\tau=\frac{\sqrt{\theta^2-1}}{\theta-1}=\sigma+\sqrt{\sigma^2\!-\!1},&&
\hspace{37pt}\eta=-\frac{\sqrt{\psi^2-1}}{\psi-1}=\frac{\tau\,\xi-1}{\xi-\tau}.
\end{eqnarray}

As mentioned in Section \ref{sec:prelim}, the quadratic transformation
identifies 4 branches on the $(\frac{A}2,\frac{A}2,\frac{B}2,\frac{B}2+1)$
side with 2 branches on the $(0,0,B,A+1)$ side. The branch permutations
on the $(\frac{A}2,\frac{A}2,\frac{B}2,\frac{B}2+1)$ side are realized by the
following transformations:
\begin{eqnarray*}
\hspace{-9pt}&\bullet&\tau\mapsto -\tau,\quad \sigma\mapsto-\sigma,\quad
\xi\mapsto\frac{\sigma\xi-1}{\sigma-\xi},\quad
 \mbox{conjugate $\sqrt{\sigma^2-1}$, $\sqrt{\theta^2-1}$};\\
\hspace{-9pt}&\bullet&\eta\mapsto -\eta,\quad
\xi\mapsto\frac{\sigma\xi-1}{\xi-\sigma},\quad
 \mbox{conjugate $\sqrt{\psi^2-1}$};\\
\hspace{-9pt}&\bullet&\tau\mapsto -\tau,\quad \eta\mapsto -\eta, \quad
\sigma\mapsto-\sigma,\quad \xi\mapsto-\xi,\quad
 \mbox{conjugate all square roots.}
\end{eqnarray*}
The branch permutation on the $(0,0,B,A+1)$ side is realized by
\begin{eqnarray*}
&\bullet&\tau\mapsto1/\tau,\quad \eta\mapsto1/\eta,\quad
\theta\mapsto-\theta,\quad \psi\mapsto-\psi,\quad \mbox{conjugate all square
roots.}
\end{eqnarray*}
Interchanging $A$ and $B$ is realized by:
\begin{eqnarray*}
\quad\,&\bullet&\eta\mapsto\frac{\tau}{\eta},\quad
\psi\mapsto\frac{\theta\psi-1}{\psi-\theta},\quad \sqrt{\psi^2-1}\mapsto
\frac{\sqrt{\theta^2-1}\sqrt{\psi^2-1}}{\psi-\theta},\quad \xi\mapsto
\frac{2\sigma-\xi+1}{\xi+1}.
\end{eqnarray*}

Now we assume the setting of Lemma \ref{lem:main}.
As usual, we identify $A=a+b-1$, $B=b-a$.
The functions $G_1,G_3$ are related by the interchange of $A$ and $B$.
In particular:
\begin{equation}
T_2=\frac12+\frac{\theta}{2\,\sqrt{\theta^2-1}},\qquad\quad
G_1(T_2)=\frac12+\frac{\theta-\psi+\sqrt{\psi^2-1}}{2\,\sqrt{\theta^2-1}}.
\end{equation}
Let us define $\varphi$ by $g_0=(1+\varphi)/2$. Then
\begin{eqnarray}
g_2(t_2)\equal\frac12-\frac{\psi\varphi^2-2\varphi+\psi}
{2\left(\varphi^2-2\psi\varphi+1\right)},\\
G_2(T_2)\equal\frac12+\frac{\theta-\psi+\sqrt{\psi^2-1}}{2\sqrt{\theta^2-1}}
+\frac{(\varphi-\psi)\left(\psi-\sqrt{\psi^2-1}\right)}
{\left(\varphi-\psi+\sqrt{\psi^2-1}\right)\sqrt{\theta^2-1}}.
\end{eqnarray}
Formula (\ref{eq2:mainsol}) can be rewritten as follows:
\begin{eqnarray}
\label{eq1:mainsol}
G_0(T_2)\equal\frac12+\frac{\theta-\psi+\sqrt{\psi^2-1}}{2\,\sqrt{\theta^2-1}}
+\frac{(a-b+1)(\varphi-\psi)\left(\psi-\sqrt{\psi^2-1}\right)}
{2\left(a\sqrt{\psi^2-1}-(b-1)(\varphi-\psi)\right)\sqrt{\theta^2-1}}.
\end{eqnarray}
Even more concisely, we can write
\begin{equation} \label{eq1:maint2}
G_0(T_2)=\frac12+\frac{\theta}{2\,\sqrt{\theta^2-1}}+
\frac{a(\psi\phi-1-\phi\sqrt{\psi^2-1})}
{2\left(a\sqrt{\psi^2-1}-(b-1)(\varphi-\psi)\right)\sqrt{\theta^2-1}}.
\end{equation}

\end{document}